\documentclass[a4paper]{siamltex}

\usepackage[english]{babel}
\usepackage{amssymb,amsfonts,amsmath,latexsym}
\usepackage[dvips]{graphicx}
\usepackage{algorithm}
\usepackage{algorithmic}
\usepackage{color}
\usepackage{url}
\usepackage{hyperref}
\newcommand{\bm}[1]{\mathbf{#1}}
\newcommand{\R}{{\mathbb{R}}}
\newcommand{\cN}{{\mathcal{N}}}
\newcommand{\cP}{{\mathcal{P}}}
\newcommand{\cQ}{{\mathcal{Q}}}
\newcommand{\cI}{{\mathcal{I}}}
\DeclareMathOperator{\Span}{span}

\newtheorem{remark}[theorem]{Remark}
\newtheorem{example}[theorem]{Example}

\begin{document}

\title{A doubly relaxed minimal-norm Gauss--Newton method for underdetermined 
nonlinear least-squares problems
%placed here. General acknowledgments should be placed at the end of the 
%article.}
%\footnotetext[0]{Version \today}
}

\author{Federica Pes\thanks{Department of Mathematics and Computer Science, via 
Ospedale 72, 09124 Cagliari, Italy, \texttt{federica.pes@unica.it, 
rodriguez@unica.it}}
\and
Giuseppe Rodriguez\footnotemark[1]
}

\maketitle

\begin{abstract}
When a physical system is modeled by a nonlinear function, the unknown 
parameters can be estimated by fitting experimental observations by a
least-squares approach.
Newton's method and its variants are often used to solve problems of this type.
In this paper, we are concerned with the computation of the minimal-norm
solution of an underdetermined nonlinear least-squares problem.
We present a Gauss--Newton type method, which relies on
two relaxation parameters to ensure convergence, and which incorporates a
procedure to dynamically estimate the two parameters, as 
well as the rank of the Jacobian matrix, along the iterations.
Numerical results are presented.
\end{abstract}

\begin{keywords}
nonlinear least-squares problem, minimal-norm solution, Gauss--Newton method, 
parameter estimation
\end{keywords}

\begin{AMS}
65H10, 65F22
\end{AMS}

\section{Introduction}\label{intro}

Let us assume that $F(\bm{x})=[F_1(\bm{x}),\ldots,F_m(\bm{x})]^T$
is a nonlinear twice continuously Frech\'et-differentiable
function with values in $\R^m$, for any $\bm{x}\in\R^n$.
For a given $\bm{b}\in\R^m$, we consider the nonlinear least-squares data
fitting problem
\begin{equation}
\min_{\bm{x}\in\R^n} \|\bm{r}(\bm{x})\|^2, \qquad \bm{r}(\bm{x}) = 
F(\bm{x})-\bm{b},
\label{nonlinL2}
\end{equation}
where $\|\cdot\|$ denotes the Euclidean norm and $\bm{r}(\bm{x})=
\left[r_1(\bm{x}),\ldots,r_m(\bm{x})\right]^T$ is the residual
vector function between the model expectation $F(\bm{x})$ and the vector
$\bm{b}$ of measured data.
The solution to the nonlinear least-squares problem gives the best model fit to
the data in the sense of the minimum sum of squared errors.
A common choice for solving a nonlinear least-squares
problem consists of applying Newton's method and its variants, such as the
Gauss--Newton method~\cite{bjo96,hansen2012least,ortega1970}.

\medskip

The Gauss--Newton method is based on the construction of a sequence of linear
approximations to $\bm{r}(\bm{x})$. Chosen an initial point $\bm{x}^{(0)}$ and
denoting by $\bm{x}^{(k)}$ the current approximation, then the new approximation
is 
\begin{equation}
\label{iter}
\bm{x}^{(k+1)}=\bm{x}^{(k)}+\bm{s}^{(k)}, \qquad k=0,1,2,\ldots,
\end{equation}
where the step $\bm{s}^{(k)}$ is computed as a solution to the linear
least-squares problem
\begin{equation}
\label{gauss-newt}
\min_{\bm{s}\in\R^n} \| J(\bm{x}^{(k)})\bm{s}+\bm{r}(\bm{x}^{(k)}) \|^2.
\end{equation}
Here $J(\bm{x})$ represents the Jacobian matrix of the function $F(\bm{x})$.

The solution to~\eqref{gauss-newt} may not be unique: this happens 
when the matrix $J(\bm{x}^{(k)})$ does not have full column rank, in
particular, when $m<n$.
To make the solution unique, the new iterate $\bm{x}^{(k+1)}$ is often obtained
by solving the following minimal-norm linear least-squares problem
\begin{equation}
\begin{cases}
\displaystyle\min_{\bm{s}\in\R^n}\|\bm{s}\|^2 \\
\displaystyle \bm{s} \in \bigl\{ \arg \min_{\bm{s}\in\R^n}
\|J(\bm{x}^{(k)})\bm{s}+\bm{r}(\bm{x}^{(k)})\|^2 \bigr\},
\end{cases}
\label{mnls}
\end{equation}
where the set in the lower line contains all the solutions to problem
\eqref{gauss-newt}.

In order to select solutions exhibiting different degrees of regularity, 
the term $\|\bm{s}\|^2$ in~\eqref{mnls} is sometimes substituted by the
seminorm $\|L\bm{s}\|^2$, where $L\in \R^{p\times n}$ $(p\leq n)$ is a matrix
which incorporates available a priori information on the solution.
The case $p>n$ can be easily reduced to the previous assumption by performing
a compact $L=QR$ factorization, and substituting $L$ by the triangular 
matrix $R$.
Typically, $L$ is a diagonal weighting matrix or a discrete
approximation of a derivative operator.
For example, the matrices
\begin{equation}\label{d1d2}
D_1= \begin{bmatrix}
1 & -1 & & & \\
  & \ddots & \ddots & \\
  & & 1 & -1
\end{bmatrix} 
\quad
\text{ and } 
\quad
D_2=\begin{bmatrix}
1 & -2 & 1 & & & \\
  & \ddots & \ddots & \ddots &\\
  & & 1 & -2 & 1
\end{bmatrix},
\end{equation}
of size $(n-1)\times n$ and $(n-2)\times n$, respectively,
are approximations to the first and second derivative operators.
When a regularization matrix is introduced,
problem~\eqref{mnls} becomes 
\begin{equation}
\begin{cases}
\displaystyle \min_{\bm{s}\in\R^n}\|L\bm{s}\|^2 \\
\displaystyle \bm{s} \in \bigl\{ \arg \min_{\bm{s}\in\R^n} 
\|J(\bm{x}^{(k)})\bm{s}+\bm{r}(\bm{x}^{(k)})\|^2 \bigr\}.
\end{cases}
\label{Lmnls}
\end{equation}

Both~\eqref{mnls} and~\eqref{Lmnls} impose some
kind of regularity on the update vector $\bm{s}$ for the solution 
$\bm{x}^{(k)}$ and not on the solution itself. 
The problem of imposing a regularity constraint directly on the solution
$\bm{x}$ of problem~\eqref{nonlinL2}, i.e.,
\begin{equation}\label{minnorm}
\begin{cases}
\displaystyle\min_{\bm{x}\in\R^n}\|\bm{x}\|^2 \\
\displaystyle \bm{x} \in \bigl\{ \arg \min_{\bm{x}\in\R^n}
\|F(\bm{x})-\bm{b}\|^2 \bigr\},
\end{cases}
\end{equation}
is studied
in~\cite{Eriksson96optimization,ErikssonPaperII,ErikssonReg05,pr20}.
These papers are based on the application of the damped Gauss--Newton method to
the solution of \eqref{minnorm}.
To ensure the computation of the minimal-norm solution, at the $k$th
iteration, the Gauss--Newton approximation is orthogonally projected onto the
null space of the Jacobian $J(\bm{x}^{(k)})$.
In~\cite{pr20}, the damping parameter is estimated by the Armijo--Goldstein
principle; we refer to this method as the MNGN algorithm.
In the same paper, this approach is applied to the minimization of a suitable
seminorm, and different regularization techniques are considered under the
assumption that the nonlinear function $F$ is ill-conditioned.

Unfortunately, the algorithms developed in the above papers
occasionally lack to converge.
They take the form
$$
\bm{x}^{(k+1)} = \bm{x}^{(k)} + \alpha_k \widetilde{\bm{s}}^{(k)} 
- \cP_{\cN(J_k)} \bm{x}^{(k)},
$$
where $\widetilde{\bm{s}}^{(k)}$ is the solution of \eqref{mnls}, $\alpha_k$ is
a step length, and $\cP_{\cN(J_k)}$ is the orthogonal projector onto the null
space of $J_k=J(\bm{x}^{(k)})$.
One reason for the nonconvergence of such methods is that the projection step
may cause the residual to increase considerably at particular iterations.
Moreover, the rank of $J(\bm{x}^{(k)})$ may vary as the iteration progresses,
and its incorrect estimation often leads to the presence of small singular
values for the Jacobian, which amplify computational errors.

This problem of nonconvergence is dealt with in~\cite{campbell}, 
%a paper we were not aware of when we wrote~\cite{pr20}, 
by a method which will be denoted CKB in the following.
The authors consider a convex combination of the Gauss--Newton approximation
and its orthogonal projection, and apply a relaxation parameter $\gamma_k$ to
this search direction, chosen according to a given rule.
After some manipulation, the method can be written as
\begin{equation}\label{ckb}
\bm{x}^{(k+1)} = \bm{x}^{(k)} + \widetilde{\bm{s}}^{(k)} 
- \gamma_k \cP_{\cN(J_k)} \bm{x}^{(k)}.
\end{equation}
This approach makes the computation of the minimal-norm solution more robust,
but it may not converge in some situation; see Section~\ref{steplen}.
Moreover, both the MNGN and the CKB methods suffer from serious convergence
problems caused by the variation of the rank of the Jacobian along the
iterations. The rank often drops to a small value in a neighborhood of the
solution, while the two methods consider a fixed rank, generally assumed to be
the smaller dimension of the Jacobian.

In this paper, we aim at improving the convergence of the methods presented
in~\cite{campbell} and~\cite{pr20}.
We do this by first introducing in the MNGN method a technique to estimate the
rank of the matrix $J(\bm{x}^{(k)})$ at each iteration.
This procedure has the effect of improving the convergence of the method,
reducing the possibility that the iteration diverges because of error
amplification.
Then, we introduce a second relaxation parameter for the projection term, as
well as a strategy to automatically tune it, besides the usual damping
parameter for the Gauss--Newton search direction.
This approach produces, on the average, solutions closer to optimality, i.e.,
with smaller norms, than those computed by the CKB method.
Furthermore, we consider a model profile $\overline{\bm{x}}$ for the
solution, which is useful in applications where sufficient a priori information
on the physical system under investigation is available.

The paper is structured as follows. 
In Section~\ref{n_mns}, we revise the MNGN method and reformulate Theorem 3.1
from~\cite{pr20} by introducing a model profile for the solution.
Then, we give a theoretical justification for the fact that the convergence of
the method may not be ensured.
Section~\ref{rankjac} explains how to estimate the numerical rank of the
Jacobian $J(\bm{x}^{(k)})$ at each iteration.
In Section~\ref{steplen}, we describe an algorithm which introduces a second
parameter to control the size of the correction vector that provides the
minimal-norm solution, and which estimates automatically such parameter.
In Section~\ref{n_mLns}, we extend the discussion to the minimal-$L$-norm
solution, where $L$ is a regularization matrix.
Numerical examples can be found in Section~\ref{examples}.

\section{Nonlinear minimal-norm solution}\label{n_mns}

We begin by recalling the definition of the singular value decomposition (SVD) 
of a matrix $J\in\R^{m\times n}$~\cite{gvl96}, which will be needed later.
The SVD is a matrix decomposition of the form
$$
J=U\Sigma V^T,
$$
where $U=[\bm{u}_1,\dots,\bm{u}_m]\in \R^{m\times m}$ and
$V=[\bm{v}_1,\dots,\bm{v}_n]\in \R^{n\times n}$ are matrices with orthonormal
columns and $\Sigma_{i,j}=0$ for $i\neq j$.
The nonzero diagonal elements of the matrix
$\Sigma \in \R^{m\times n}$ are the \emph{singular values}
$\sigma_1\geq\sigma_2\geq\cdots\geq\sigma_r>0$, with
$r=\rank(J)\leq\min(m,n)$.
Let $\cN(J)$ denote the null space of the matrix $J$.
It is well-known that
$$
\cN(J) := \left\{ \bm{s}\in\R^n : J\bm{s}=0 \right\}
= \Span\{ \bm{v}_{r+1},\ldots,\bm{v}_n \}.
$$

Let us now briefly review the computation of the minimal-norm solution to the
nonlinear problem~\eqref{nonlinL2} by the \emph{minimal-norm Gauss--Newton}
(MNGN) method, presented in~\cite{pr20}.
Our aim is showing the reason for the possible lack of convergence of such
method.
%A similar analysis may be developed for the CKB method from~\cite{campbell}.
Here, we extend the discussion from~\cite{pr20} by introducing a model profile
$\overline{\bm{x}}\in\R^n$, which represents an a priori estimate of the
desired solution, and formulate the problem in the form
\begin{equation}
\begin{cases}
\displaystyle\min_{\bm{x}\in\R^n}\|\bm{x}-\overline{\bm{x}}\|^2 \\
\displaystyle \bm{x} \in \bigl\{ \arg \min_{\bm{x}\in\R^n}
\|F(\bm{x})-\bm{b}\|^2 \bigr\}.
\end{cases}
\label{nonlinmnls}
\end{equation}
We consider an iterative method of the type~\eqref{iter} based on the following
first-order linearization of the problem
\begin{equation}
\begin{cases}
\displaystyle\min_{\bm{s}\in\R^n}\|\bm{x}^{(k)}-\overline{\bm{x}}+\alpha_k\bm{s}\|^2
 \\
\displaystyle \bm{s} \in \bigl\{ \arg \min_{\bm{s}\in\R^n} 
\|J_k\bm{s}+\bm{r}_k\|^2 \bigr\},
\end{cases}
\label{linmnal}
\end{equation}
where $J_k=J(\bm{x}^{(k)})$ is the Jacobian of $F$ in $\bm{x}^{(k)}$ and 
$\bm{r}_k=\bm{r}(\bm{x}^{(k)})$ is the residual vector.

The damping parameter $\alpha_k$ is indispensable to ensure the convergence of
the Gauss--Newton method. 
We estimate it by the Armijo--Goldstein principle~\cite{armijo,goldstein}, but 
it can be chosen by any strategy which guarantees a reduction in the norm of 
the residual.
In our case, the Armijo condition~\cite{armijo,dennis} implies
\[
f(\bm{x}^{(k)}+\alpha_k \widetilde{\bm{s}}^{(k)}) \leq f(\bm{x}^{(k)}) + 
\mu\alpha_k \nabla f(\bm{x}^{(k)})^T \widetilde{\bm{s}}^{(k)},
\]
where $\widetilde{\bm{s}}^{(k)}$ is determined by solving~\eqref{mnls} and 
$\mu$ is a constant in $(0, 1)$.
Since $f(\bm{x})=\frac{1}{2}\|\bm{r}(\bm{x})\|^2$ and 
$\nabla f(\bm{x})=J(\bm{x})^T\bm{r}(\bm{x})$, it reads
\[
\|\bm{r}(\bm{x}^{(k)}+\alpha_k \widetilde{\bm{s}}^{(k)})\|^2 
\leq \|\bm{r}_k\|^2 + 2\mu \alpha_k \bm{r}_k^T J_k \widetilde{\bm{s}}^{(k)}.
\]
Note that, as $\widetilde{\bm{s}}^{(k)}$ satisfies the normal equations
associated to problem \eqref{gauss-newt}, it holds
$J_k^T\bm{r}_k=-J_k^TJ_k\widetilde{\bm{s}}^{(k)}$, so that
$\bm{r}_k^TJ_k\widetilde{\bm{s}}^{(k)}=-\|J_k\widetilde{\bm{s}}^{(k)}\|^2$.
%$$
%\bm{r}_k^T J_k \widetilde{\bm{s}}^{(k)} = -\bm{r}_k^T J_kJ_k^\dagger \bm{r}_k
%= -\bm{r}_k^T J_kJ_k^\dagger J_kJ_k^\dagger \bm{r}_k
%= -\|J_kJ_k^\dagger \bm{r}_k\|^2
%= -\|J_k \widetilde{\bm{s}}^{(k)}\|^2,
%$$
%as $J_kJ_k^\dagger$ is the orthogonal projection onto the range of $J_k$.
The \emph{Armijo--Goldstein principle} \cite{bjo96,goldstein} sets 
$\mu=\frac{1}{4}$ and determines the
scalar $\alpha_k$ as the largest number in the sequence $2^{-i}$,
$i=0,1,\ldots,$ for which it holds
\begin{equation}\label{armgol}
\|\bm{r}_k\|^2 - \|\bm{r}(\bm{x}^{(k)}+\alpha_k \widetilde{\bm{s}}^{(k)})\|^2 
\geq \frac{1}{2} \alpha_k \|J_k \widetilde{\bm{s}}^{(k)}\|^2.
\end{equation}

The iteration resulting from the solution of \eqref{linmnal} 
is defined by the following theorem.

\medskip
\begin{theorem}\label{theo3.1}
Let $\bm{x}^{(k)}\in\R^n$ and 
let $\widetilde{\bm{x}}^{(k+1)}=\bm{x}^{(k)}+\alpha_k\widetilde{\bm{s}}^{(k)}$
be the Gauss--Newton iteration for~\eqref{nonlinL2}, 
where the step $\widetilde{\bm{s}}^{(k)}$ is determined by solving~\eqref{mnls}
and the step length $\alpha_k$ by the Armijo--Goldstein principle.
Then, the iteration $\bm{x}^{(k+1)}=\bm{x}^{(k)}+\alpha_k\bm{s}^{(k)}$ defined
by~\eqref{linmnal}
%starting from the same point $\bm{x}^{(k)}$, 
is given by
\begin{equation}\label{tesi}
\bm{x}^{(k+1)} = \widetilde{\bm{x}}^{(k+1)} - V_2V_2^T 
\bigl(\bm{x}^{(k)}-\overline{\bm{x}}\bigr), 
\end{equation}
where $\rank(J_k)=r_k$ and the columns of the matrix 
$V_2=[\bm{v}_{r_k+1},\ldots,\bm{v}_n]$
are orthonormal vectors in $\R^n$ spanning the null space of $J_k$.
\end{theorem}

\smallskip
\begin{proof}
The proof follows the pattern of that of Theorem 3.1 in~\cite{pr20}.
Let $U\Sigma V^T$ be the singular value decomposition of the matrix $J_k$.
The upper-level problem in~\eqref{linmnal} can be expressed as
\[
\|\bm{x}^{(k)}-\overline{\bm{x}}+\alpha_k\bm{s}\|^2 =
\|V^T(\bm{x}^{(k)}-\overline{\bm{x}}+\alpha_k\bm{s})\|^2 = 
\|\alpha_k\bm{y}+\bm{z}^{(k)}\|^2,
\]
with $\bm{y}=V^T\bm{s}$ and 
$\bm{z}^{(k)}=V^T\left(\bm{x}^{(k)}-\overline{\bm{x}}\right)$.
Replacing $J_k$ by its SVD and setting $\bm{g}^{(k)}=U^T\bm{r}_k$, 
we can rewrite~\eqref{linmnal} as the following diagonal linear 
least-squares problem
\[
\begin{cases}
\displaystyle \min_{\bm{y}\in\R^n}\|\alpha_k\bm{y}+\bm{z}^{(k)}\|^2 \\
\displaystyle \bm{y} \in \bigl\{ \arg \min_{\bm{y}\in\R^n}
\|\Sigma\bm{y}+\bm{g}^{(k)}\|^2 \bigr\}.
\end{cases}
\]
Solving the lower-level minimization problem uniquely determines the components
$y_i=-\sigma_i^{-1}g^{(k)}_i$, $i=1,\ldots,r_k$, while the entries $y_i$,
$i=r_k+1,\ldots,n$, are left undetermined.
Their values can be found by solving the upper-level problem.
From
\[
\|\alpha_k\bm{y}+\bm{z}^{(k)}\|^2 = \sum_{i=1}^{r_k}
\left(-\alpha_k\frac{g^{(k)}_i}{\sigma_i}+z^{(k)}_i\right)^2 + \sum_{i=r_k+1}^n
\left(\alpha_k y_i+z^{(k)}_i\right)^2,
\]
we obtain $y_i=-\frac{z^{(k)}_i}{\alpha_k} =-\frac{1}{\alpha_k} 
\bm{v}_i^T(\bm{x}^{(k)}-\overline{\bm{x}})$, $i=r_k+1,\ldots,n$.
Then, the solution to~\eqref{linmnal}, that is, the next approximation to the
solution of~\eqref{nonlinmnls}, is
$$
\bm{x}^{(k+1)} %&= \bm{x}^{(k)} + \alpha_k\bm{s}
= \bm{x}^{(k)} + \alpha_k V\bm{y} 
= \bm{x}^{(k)} - \alpha_k\sum_{i=1}^{r_k} \frac{g^{(k)}_i}{\sigma_i} \bm{v}_i
-\sum_{i=r_k+1}^n (\bm{v}_i^T(\bm{x}^{(k)}-\overline{\bm{x}})) \bm{v}_i,
$$
where the last summation can be written in matrix form as
$V_2V_2^T\left(\bm{x}^{(k)}-\overline{\bm{x}}\right)$, and the columns of
$V_2=[\bm{v}_{r_k+1},\ldots,\bm{v}_n]$ are a basis for $\cN(J_k)$.

It is immediate (see \cite[Theorem~3.1]{pr20}) to prove that
\[
\widetilde{\bm{x}}^{(k+1)} = \bm{x}^{(k)} + \alpha_k\widetilde{\bm{s}}^{(k)} 
= \bm{x}^{(k)} - \alpha_k\sum_{i=1}^{r_k} \frac{g^{(k)}_i}{\sigma_i} \bm{v}_i,
\]
from which \eqref{tesi} follows.
\end{proof}
\medskip

%\medskip
%\begin{corollary}
%If $\bm{x}^{(k)}$ is orthogonal to the null space of $J_k$, then the solution 
%$\bm{x}^{(k+1)}$ of~\eqref{linmnal} is the same as that of~\eqref{mnls}. 
%\end{corollary}
%\medskip

Summarizing, the MNGN method consists of the iteration
\[
\bm{x}^{(k+1)}=\bm{x}^{(k)}+\alpha_k\bm{s}^{(k)},
\]
where the step is
$$
\bm{s}^{(k)}= \widetilde{\bm{s}}^{(k)} - \frac{1}{\alpha_k} \bm{t}^{(k)},
$$
with
\begin{equation} \label{stepdef}
\widetilde{\bm{s}}^{(k)} = -\sum_{i=1}^{r_k} \frac{g^{(k)}_i}{\sigma_i} 
\bm{v}_i,
\qquad
\bm{t}^{(k)} = V_2V_2^T\bigl(\bm{x}^{(k)} -\overline{\bm{x}}\bigr).
\end{equation}
Since $\cP_{\cN(J_k)}=V_2V_2^T$ is the orthogonal projector onto $\cN(J_k)$, 
the above theorem states that the $(k+1)$th iterate of the MNGN method is 
orthogonal to the null space of $J_k$.
\smallskip

Theorem~\ref{theo3.1} shows that the correction vector $\bm{t}^{(k)}$ defined
in~\eqref{stepdef}, which allows to compute the minimal-norm solution at each
step, is not damped by the parameter $\alpha_k$.
As a result, in some numerical examples, the method fails to converge
because projecting the solution orthogonally to the null space of $J_k$ causes
the residual to increase.
To understand how this can happen, a second-order analysis of the objective
function is required.

The second-order Taylor 
approximation to the function $f(\bm{x})=\frac{1}{2}\|\bm{r}(\bm{x})\|^2$ at
$\bm{x}^{(k+1)}=\bm{x}^{(k)}+\alpha\bm{s}$ is
\begin{equation}\label{second}
f(\bm{x}^{(k+1)}) \simeq f(\bm{x}^{(k)})+ \alpha\nabla f(\bm{x}^{(k)})^T 
\bm{s} + \frac{1}{2} \alpha^2 \bm{s}^T \nabla^2 f(\bm{x}^{(k)}) \bm{s}.
\end{equation}
The gradient and the Hessian of $f(\bm{x})$, written in matrix form, are given
by
\[
\nabla f(\bm{x})= J(\bm{x})^T \bm{r}(\bm{x}), \qquad
\nabla^2 f(\bm{x})= J(\bm{x})^T J(\bm{x}) + \cQ(\bm{x}), 
\]
where
\[
\cQ(\bm{x})= \sum_{i=1}^m r_i(\bm{x})\nabla^2 r_i(\bm{x}),
\]
and $\nabla^2 r_i(\bm{x})$ is the Hessian matrix of $r_i(\bm{x})$.
By replacing the expression of $f$ and 
$\alpha\bm{s}=\alpha\widetilde{\bm{s}}-\bm{t}$
in~\eqref{second}, where $\widetilde{\bm{s}}$ is the Gauss--Newton step and
$\bm{t}$ is in the null space of $J_k$, 
and letting $\cQ_k=\cQ(\bm{x}^{(k)})$,
the following approximation is obtained
\[
\begin{aligned}
\frac{1}{2}\|\bm{r}_{k+1}\|^2 
&\simeq \frac{1}{2}\|\bm{r}_k\|^2 + \alpha \bm{r}_k^T J_k 
\bm{s} + \frac{1}{2} \alpha^2 \bm{s}^T \left(J_k^T J_k + \cQ_k\right)\bm{s} \\
&= \frac{1}{2}\|\bm{r}_k\|^2 + \alpha \bm{r}_k^T J_k 
\widetilde{\bm{s}} + \frac{1}{2} \alpha^2 \widetilde{\bm{s}}^T \left(J_k^T J_k 
+ \cQ_k\right)\widetilde{\bm{s}} - \alpha\bm{t}^T \cQ_k \widetilde{\bm{s}} + 
\frac{1}{2}\bm{t}^T \cQ_k \bm{t}.
\end{aligned}
\]

The first two terms containing second derivatives (the matrix $\cQ_k$) are
damped by the $\alpha$ parameter.
If the function $F$ is mildly nonlinear, the third term $\frac{1}{2}\bm{t}^T
\cQ_k \bm{t}$ is negligible. In the presence of a strong nonlinearity, its
contribution to the residual is significant and may lead to its growth. This
shows that a damping parameter is required to control the step length for both
the Gauss--Newton step $\widetilde{\bm{s}}$ and the correction vector $\bm{t}$.
If a relaxation parameter is introduced for $\bm{t}$, Theorem~\ref{theo3.1}
implies that the minimal-norm solution of \eqref{linmnal} can only be
approximated.

\begin{remark}\label{smallex}\rm
We report a simple low dimensional example for which the MNGN method may not
converge. Let us consider the function $F:\R^2 \rightarrow \R$ defined by
$$
F(\bm{x})=\delta^2 \left[ (x_1-\gamma)^2+(x_2-\gamma)^2 \right]-1,
$$
depending on the parameters $\delta,\gamma\in\R$.
Since the Hessian matrix of the residual is given by
$$
\nabla^2 r(\bm{x})=
\begin{bmatrix} 2\delta^2 & 0\\ 0 & 2\delta^2 \end{bmatrix},
$$
the second-order term $\frac{1}{2}\bm{t}^T \cQ_k \bm{t}$ is not negligible, in
general, when $\delta$ is relatively large.
For example, setting $\delta=0.7$, $\gamma=2$, and choosing an initial vector
$\bm{x}^{(0)}$ with random components in $(-5,5)$, the MNGN method converges 
with a large number of the iterations (350 on average).
Setting $\delta=0.75$, the same method does not converge within 500 iterations.
\end{remark}

\section{Estimating the rank of the Jacobian}\label{rankjac}

In order to apply Theorem~\ref{theo3.1} to computing the minimal-norm solution
by \eqref{tesi}, the rank of the Jacobian matrix $J_k=J(\bm{x}^{(k)})$ should
be known in advance.
As the rank may vary along the iterations, we set $r_k=\rank(J_k)$.
The knowledge of $r_k$ for each $k=0,1,\ldots$, is not generally available,
making it necessary to estimate its value at each iteration step, to avoid
nonconvergence or a breakdown of the algorithm.

In such situations, it is common to consider the numerical rank $r_{\epsilon,k}$
of $J_k$, sometimes denoted as $\epsilon$-rank, where $\epsilon$ represents a
chosen tolerance. 
The numerical rank is defined in terms of the singular values $\sigma_i^{(k)}$
of $J_k$, as the integer $r_{\epsilon,k}$ such that
\[
\sigma_{r_{\epsilon,k}}^{(k)}>\epsilon\geq \sigma_{r_{\epsilon,k}+1}^{(k)}.
\]
%We associate to the numerical rank the \emph{numerical null space} of $J_k$,
%defined as
%\[
%\cN_{r_{\epsilon,k}}(J_k)=\Span\left\{\bm{v}_{r_{\epsilon,k}+1}, \ldots,
%\bm{v}_n \right\}.
%\]
Theorem~\ref{theo3.1} can be adapted to this setting, by simply replacing at
each iteration the rank $r_k$ with the numerical rank $r_{\epsilon,k}$.

Determining the numerical rank is a difficult task for discrete ill-posed
problems, in which the singular values decay monotonically to zero.
In such a case, the
numerical rank plays the role of a regularization parameter and is estimated by
suitable methods, which often require information about the noise level and
type; see, e.g., \cite{Hansen,rr13}.

When the problem is locally rank-deficient, meaning that the rank of
$J(\bm{x})$ depends on the evaluation vector $\bm{x}$, the numerical rank
$r_{\epsilon,k}$ can be determined, in principle, by choosing a suitable value
of $\epsilon$.
Numerical experiments show that a fixed value of $\epsilon$ does not always
lead to a correct estimation of $r_{\epsilon,k}$, and that it is preferable to
determine the $\epsilon$-rank by searching for a sensible gap between
$\sigma_{r_{\epsilon,k}}^{(k)}$ and $\sigma_{r_{\epsilon,k}+1}^{(k)}$.

To locate such a gap, we adopt a heuristic approach already applied
in~\cite{cnrr20} for the same purpose, in a different setting.
At each step, we compute the ratios
\[
\rho_i^{(k)} = \frac{\sigma_i^{(k)}}{\sigma_{i+1}^{(k)}}, 
\qquad i=1,2,\ldots,q-1,
\]
where $q=\min(m,n)$.
Then, we consider the index set
\[
\cI_k = \left\{ i\in\{1,2,\ldots,q-1\} : \rho_i^{(k)}>R \text{ and } 
\sigma_i^{(k)}>\tau \right\}.
\]
An index $i$ belongs to $\cI_k$ if there is a significant ``jump'' between
$\sigma_i^{(k)}$ and $\sigma_{i+1}^{(k)}$, and $\sigma_i^{(k)}$ is
numerically nonzero.
If the set $\cI_k$ is empty, we set $r_{\epsilon,k}=q$.
Otherwise, we consider
\begin{equation}\label{rankest}
\rho_j^{(k)}=\max_{i\in\cI_k}\rho_i^{(k)},
\end{equation}
and we define $r_{\epsilon,k}=j$.
This amounts to selecting the largest gap between ``large'' and ``small'' 
singular values.
In our numerical simulations, we set $R=10^2$ and $\tau=10^{-8}$.
We observed that the value of these parameters is not critical for problems
characterized by a rank deficient Jacobian.
Estimating the rank becomes increasingly difficult as the gap between 
``large'' and ``small'' singular values gets smaller.
This condition usually corresponds to ill-conditioned problems, which require
specific regularization methods.

\section{Choosing the projection step length}\label{steplen}

The occasional nonconvergence in the computation of the minimal-norm solution
to a nonlinear least-squares problem was discussed in~\cite{campbell}, where
the authors propose an iterative method based on a convex combination of the
Gauss--Newton and the minimal-norm Gauss--Newton iterates, which we denote by
CKB.
Following our notation, it can be expressed in the form
\begin{equation}
\label{camp}
\bm{x}^{(k+1)} = \left(1-\gamma_k\right)\left[\bm{x}^{(k)} +
\widetilde{\bm{s}}^{(k)}\right] + \gamma_k\left[\bm{x}^{(k)} +
\widetilde{\bm{s}}^{(k)} - V_2V_2^T\bm{x}^{(k)}\right],
\end{equation}
where the parameters $\gamma_k\in[0,1]$, for $k=0,1,\ldots$, form a sequence
converging to zero.
The standard Gauss--Newton method is obtained by setting $\gamma_k=0$, while
$\gamma_k=1$ leads to the minimal-norm Gauss--Newton method.
In their numerical examples, the authors adopt the sequences
$\gamma_k=(0.5)^{k+1}$ and $\gamma_k=(0.5)^{2^k}$.

It is immediate to rewrite~\eqref{camp} in the form \eqref{ckb},
showing that the method proposed in~\cite{campbell} is equivalent to the
application of the undamped Gauss--Newton method, whose convergence is not
theoretically guaranteed~\cite{bjo96}, with a damped correction to favor the
decrease of the norm of the solution. 
The numerical experiments reported in the paper show that the minimization of
the residual is sped up if $\gamma_k$ quickly converges to zero, while the norm
of the solution decreases faster if $\gamma_k$ has a slower decay. 
The choice of the sequence of parameters appears to be critical to tune the
performance of the algorithm, and no adaptive choice for $\gamma_k$ is
proposed.
\smallskip

In this paper, we propose to introduce a second relaxation parameter,
$\beta_k$, to control the step length of the minimal-norm correction
$\bm{t}^{(k)}$ defined in~\eqref{stepdef}.
The new iterative method is denoted by MNGN2 and it takes the form
\begin{equation}\label{mngn2}
\bm{x}^{(k+1)} = \bm{x}^{(k)} + \alpha_k \widetilde{\bm{s}}^{(k)} 
- \beta_k \bm{t}^{(k)},
\end{equation}
where $\widetilde{\bm{s}}^{(k)}$ is the step vector produced by the
Gauss--Newton method and $\bm{t}^{(k)}$ is the projection vector which makes the
norm of $\bm{x}^{(k+1)}$ minimal, without changing the value of the linearized
residual.

The second-order analysis reported at the end of Section~\ref{n_mns} may be
adapted for the CKB method \eqref{ckb}.
It shows that neither the CKB nor the MNGN method are guaranteed to converge,
as both the Gauss--Newton search direction and the projection step should be
damped to ensure that the residual decreases.
The MNGN2 method locally converges if $\alpha_k$ and $\beta_k$ are suitably
chosen, but it will recover the minimal-norm solution only if 
$\beta_k\simeq 1$ for $k$ close to convergence.

Our numerical tests showed that it is important to choose both $\alpha_k$ and
$\beta_k$ adaptively along the iterations.
A simple solution is to let $\beta_k=\alpha_k$ and estimate $\alpha_k$ by the
Armijo--Goldstein principle~\eqref{armgol}, with
$\bm{s}^{(k)}=\widetilde{\bm{s}}^{(k)}-\bm{t}^{(k)}$ in place of
$\widetilde{\bm{s}}^{(k)}$.
This approach proves to be effective in the computation of the minimal-norm
solution, but its convergence is often rather slow.
To speed up iteration we propose a procedure to adaptively choose the value of
$\beta_k$.

\begin{algorithm}
\caption{Outline of the MNGN2 method.}
\label{algobeta}
\begin{algorithmic}[1]
\REQUIRE nonlinear function $F$, data vector $\bm{b}$, 
\REQUIRE initial solution $\bm{x}^{(0)}$, model profile $\overline{\bm{x}}$,
tolerance $\eta$ for residual increase 
\ENSURE approximation $\bm{x}^{(k+1)}$ of minimal-norm least-squares solution
\STATE $k=0$, $\beta=1$
\REPEAT
	\STATE $k=k+1$
	\STATE estimate $r_k=\rank(J(\bm{x}^{(k)}))$ by \eqref{rankest}
	\STATE compute $\widetilde{\bm{s}}^{(k)}$ by the Gauss--Newton method 
	\eqref{gauss-newt}
	\STATE compute $\alpha_k$ by the Armijo--Goldstein principle \eqref{armgol}
	\STATE compute $\bm{t}^{(k)}$ by \eqref{stepdef}
	\IF {$\beta<1$}
		\STATE $\beta=2\beta$ \label{beta2}
	\ENDIF
	\STATE 
	$\widetilde{\bm{x}}^{(k+1)}=\bm{x}^{(k)}+\alpha_k\widetilde{\bm{s}}^{(k)}$ 
	\STATE 
	$\widetilde{\rho}_{k+1}=\|F(\widetilde{\bm{x}}^{(k+1)})-\bm{b}\|+\varepsilon_M$
	 \label{linevareps}
	\STATE $\bm{x}^{(k+1)}=\widetilde{\bm{x}}^{(k+1)}-\beta\bm{t}^{(k)}$
	\STATE $\rho_{k+1}=\|F(\bm{x}^{(k+1)})-\bm{b}\|$ 
	\WHILE {$(\rho_{k+1}> 
	\widetilde{\rho}_{k+1}+\delta(\widetilde{\rho}_{k+1},\eta))$
	\AND ($\beta>10^{-8}$)}
		\STATE $\beta=\beta/2$
		\STATE $\bm{x}^{(k+1)}=\widetilde{\bm{x}}^{(k+1)}-\beta\bm{t}^{(k)}$
		\STATE $\rho_{k+1}=\|F(\bm{x}^{(k+1)})-\bm{b}\|$ 
	\ENDWHILE
	\STATE $\beta_k=\beta$
\UNTIL {convergence}
\end{algorithmic}
\end{algorithm}

This procedure is outlined in Algorithm~\ref{algobeta}.
Initially, we set $\beta=1$.
At each iteration, we compute the residual at the Gauss--Newton iteration
$\widetilde{\bm{x}}^{(k+1)}$ and at the tentative iteration
$\bm{x}^{(k+1)}=\widetilde{\bm{x}}^{(k+1)}-\beta\bm{t}^{(k)}$.
Subtracting the vector $\beta\bm{t}^{(k)}$ may cause the residual to increase.
We accept such an increase if
\begin{equation}
\label{condiz}
\|\bm{r}(\bm{x}^{(k+1)})\| \leq \|\bm{r}(\widetilde{\bm{x}}^{(k+1)})\| +
\delta\bigl(\|\bm{r}(\widetilde{\bm{x}}^{(k+1)})\|,\eta\bigr),
\end{equation}
where $\delta(\rho,\eta)$ is a function determining the maximal increase
allowed in the residual $\rho=\|\bm{r}(\widetilde{\bm{x}}^{(k+1)})\|$, and
$\eta>0$ is a chosen tolerance.
On the contrary, $\beta$ is halved and the residual is recomputed 
until~\eqref{condiz} is verified or $\beta$ becomes excessively small.
To allow $\beta$ to increase, we tentatively double it at each iteration (see
line~\ref{beta2} in the algorithm) before applying the above procedure.
At line~\ref{linevareps} of the algorithm we add the machine epsilon
$\varepsilon_M$ to the actual residual $\widetilde{\rho}_{k+1}$ to avoid that
$\delta(\widetilde{\rho}_{k+1},\eta)$ becomes zero.

A possible choice for the value of the residual increase is
$\delta(\rho,\eta)=\eta\rho$, with $\eta$ suitably chosen.
Our experiments showed that it is possible to find, by chance, a value of
$\eta$ which produces good results, but its choice is strongly dependent on the
particular example.
We also noticed that, in cases where the residual stagnates, accepting a
large increase in the residual may lead to nonconvergence.
In such situations, a fixed multiple of the residual is not well suited to
model its increase.
Indeed, if the residual is large, one is prone to accept only a small increase,
while if the residual is very small, a relatively large growth may be
acceptable.

To overcome these difficulties, we consider $\delta(\rho,\eta)=\rho^\eta$, and
choose $\eta$ at each step by the adaptive procedure described in
Algorithm~\ref{algoeta}.
When at least $k_{\text{res}}$ iterations have been performed, we compute the
linear polynomial which fits the logarithm of the last $k_{\text{res}}$
residuals in the least-squares sense.
To detect if the residual stagnates or increases, we check
if the slope $M$ of the regression line exceeds $-10^{-2}$.
If this happens, the value of $\eta$ is doubled.
The effect on the algorithm is to enhance the importance of the decrease of the
residual and reduce that of the norm.
To recover a sensible decrease in the norm,
if at a subsequent step the residual reduction accelerates 
(e.g., $M<-\frac{1}{2}$), the value of $\eta$ is halved.
In our experiments, we initialize $\eta$ to $\frac{1}{8}$ and set
$k_{\text{res}}=5$.

\begin{remark}\label{complexity}\rm
The adaptive estimation of $\delta(\rho,\eta)$
does not significantly increase the complexity of
Algorithm~\ref{algobeta}, as line~\ref{rline} of Algorithm~\ref{algoeta}
implies the solution of a $2\times 2$ linear system whose matrix is fixed and
can be computed in advance, while forming the right-hand side requires
$4k_{\text{res}}$ floating point operations.
\end{remark}

\begin{algorithm}
\caption{Adaptive determination of the residual increase $\delta(\rho,\eta)$.}
\label{algoeta}
\begin{algorithmic}[1]
\REQUIRE actual residual $\rho=\|\bm{r}(\widetilde{\bm{x}}^{(k+1)})\|$,
starting tolerance $\eta$
\REQUIRE iteration index $k$, residuals
$\theta_j=\|\bm{r}(\widetilde{\bm{x}}^{(k-k_{\text{res}}+j)})\|$,
$j=1,\ldots,k_{\text{res}}$
\ENSURE residual increase $\delta(\rho,\eta)$
\STATE $M_{\text{min}}=-10^{-2}$, $M_{\text{max}}=-\frac{1}{2}$
\IF {$k \geq k_{\text{res}}$}
	\STATE compute regression line $p_1(t)=Mt+N$ of
		$(j,\log(\theta_j))$, $j=1,\ldots,k_{\text{res}}$ \label{rline}
	\IF {$M>M_{\text{min}}$}
		\STATE $\eta=2\eta$
	\ELSIF {$M<M_{\text{max}}$}
		\STATE $\eta=\eta/2$
	\ENDIF
\ENDIF
\STATE $\delta(\rho,\eta)=\rho^\eta$
\end{algorithmic}
\end{algorithm}

To detect convergence, we interrupt the iteration as soon as 
\begin{equation}\label{conver}
\|\bm{x}^{(k+1)}-\bm{x}^{(k)}\|<\tau\|\bm{x}^{(k+1)}\|
\qquad\text{or}\qquad
\|\alpha_k\widetilde{\bm{s}}^{(k)}\|<\tau,
\end{equation}
or when a fixed number of iteration $N_\text{max}$ is exceeded.
The second stop condition in \eqref{conver} detects the slow progress of the
relaxed Gauss--Newton iteration algorithm. This often happens close
to the solution. The stop tolerance is set to $\tau=10^{-8}$.

\section{Nonlinear minimal-$\boldsymbol{L}$-norm solution}\label{n_mLns}

The introduction of a regularization matrix $L\in\R^{p\times n}$, $p\leq n$, in
least-squares problems was originally connected to the numerical treatment of
linear discrete ill-posed problems, and in particular to Tikhonov
regularization.
The use of a regularization matrix is also justified in underdetermined
least-squares problems to select a solution with particular features, such as
smoothness or sparsity, among the infinitely many possible solutions.

While in \eqref{Lmnls} the seminorm $\|L\bm{s}\|$ is minimized over all
the updating vectors $\bm{s}$ which minimize the linearized residual, here
we seek to compute the minimal-$L$-norm solution to the nonlinear 
problem~\eqref{nonlinL2}, that is the vector $\bm{x}$ which solves
the constrained problem
\begin{equation}
\begin{cases}
\displaystyle\min_{\bm{x}\in\R^n}\|L(\bm{x}-\overline{\bm{x}})\|^2 \\
\displaystyle \bm{x} \in \bigl\{ \arg \min_{\bm{x}\in\R^n}
\|F(\bm{x})-\bm{b}\|^2 \bigr\}.
\end{cases}
\label{Lnonlinmnls}
\end{equation}
Similarly to Section~\ref{n_mns}, we consider an iterative method of the
type~\eqref{iter}, where the step $\bm{s}^{(k)}$ is the solution of the
linearized problem
\begin{equation}
\begin{cases}
\displaystyle\min_{\bm{s}\in\R^n}\|L(\bm{x}^{(k)}-\overline{\bm{x}}+\alpha\bm{s})\|^2
 \\
\displaystyle \bm{s} \in \bigl\{ \arg \min_{\bm{s}\in\R^n} 
\|J_k\bm{s}+\bm{r}_k\|^2 \bigr\}.
\end{cases}
\label{Llinmnls}
\end{equation}
We will denote the iteration resulting from the solution of \eqref{Llinmnls} as
the \emph{minimal-$L$-norm Gauss--Newton} (MLNGN) method.

We recall the definition of the generalized singular value decomposition (GSVD) 
of a matrix pair $(J,L)$ \cite{gvl96}.
Let $J\in\R^{m\times n}$ and $L\in\R^{p\times n}$ be matrices with
$\rank(J)=r$ and $\rank(L)=p$. Assume that $m+p\geq n$ and 
$$
\rank\left(\begin{bmatrix} J\\ L \end{bmatrix}\right)=n,
$$
which corresponds to requiring that $\mathcal{N}(J)\cap\mathcal{N}(L)=\{0\}$.
The GSVD of the matrix pair $(J,L)$ is defined as the factorization
$$
J=U\Sigma_J W^{-1}, \qquad
L=V\Sigma_L W^{-1},
$$
where $U\in\R^{m\times m}$ and $V\in\R^{p\times p}$ are matrices with
orthonormal columns $\bm{u}_i$ and $\bm{v}_i$, respectively,
and $W\in \R^{n\times n}$ is nonsingular. 
If $m\geq n\geq r$, the matrices $\Sigma_J \in \R^{m\times n}$ and $\Sigma_L 
\in \R^{p\times n}$ have the form
$$
\Sigma_J=\left[\begin{array}{ccc}
O_{n-r} & & \\ & C & \\ & & I_d \\ 
\hline \\
& O_{(m-n)\times n} & 
\end{array}\right], \qquad
\Sigma_L= \left[\begin{array}{cc|c}
I_{p-r+d} & & \\ & & O_{p\times d}  \\ & S &
\end{array}\right],
$$
where $d=n-p$,
\begin{equation}\label{csmat}
\begin{aligned}
C &= \diag (c_1,\ldots,c_{r-d}), \qquad & 0<c_1\leq c_2 \leq \cdots \leq 
c_{r-d} < 1, \\
S &= \diag (s_1,\ldots,s_{r-d}), \qquad & 1>s_1 \geq s_2 \geq \cdots \geq 
s_{r-d} > 0,
\end{aligned}
\end{equation}
with $c_i^2+s_i^2=1$, for $i=1,\ldots,r-d$.
The identity matrix of size $k$ is denoted by $I_k$, while $O_k$ and
$O_{k\times\ell}$ are zero matrices of size $k\times k$ and $k\times\ell$,
respectively; a matrix block has to be omitted when one of its dimensions is
zero.
The scalars $\gamma_i=\frac{c_i}{s_i}$ are called \emph{generalized singular 
values}, and they appear in nondecreasing order.

If $r\leq m<n$, the matrices $\Sigma_J \in \R^{m\times n}$ and 
$\Sigma_L \in \R^{p\times n}$ take the form
$$
\Sigma_J=\left[\begin{array}{c|ccc}
& O_{m-r} & & \\ 
O_{m\times (n-m)} & & C & \\ 
& & & I_d 
\end{array}\right], \qquad
\Sigma_L= \left[\begin{array}{cc|c}
I_{p-r+d} & & \\ & & O_{p\times d}  \\ & S &
\end{array}\right],
$$
where the blocks are defined as above.
\medskip

Let $J_k=U\Sigma_J W^{-1}$, $L=V\Sigma_L W^{-1}$ be the GSVD of the matrix pair
($J_k$,$L$). We indicate by $\bm{w}_i$ the column vectors of the matrix $W$, and
by $\widehat{\bm{w}}^j$ the rows of $W^{-1}$, that is
$$
W = [\bm{w}_1,\ldots,\bm{w}_n], \qquad
W^{-1} = \begin{bmatrix} \widehat{\bm{w}}^1 \\ \vdots \\ \widehat{\bm{w}}^n 
\end{bmatrix}.
$$
We have $\cN(J_k)=\Span(\bm{w}_1,\ldots,\bm{w}_{n-r_k})$, if 
$r_k=\rank(J_k)$; see~\cite{pr20} for a proof.

\medskip
\begin{theorem}\label{theo4.2}
Let $\bm{x}^{(k)}\in\R^n$ and let 
$\widetilde{\bm{x}}^{(k+1)}=\bm{x}^{(k)}+\alpha_k\widetilde{\bm{s}}^{(k)}$ be
the Gauss--Newton iteration for~\eqref{nonlinL2}, where the step 
$\widetilde{\bm{s}}^{(k)}$ is determined by solving~\eqref{Lmnls}
and the step length $\alpha_k$ by the Armijo--Goldstein principle.
Then, the iteration $\bm{x}^{(k+1)}=\bm{x}^{(k)}+\alpha_k\bm{s}^{(k)}$ 
for~\eqref{Llinmnls}, is given by
\begin{equation}\label{mlnsol}
\bm{x}^{(k+1)} = \widetilde{\bm{x}}^{(k+1)} - 
W_1\widehat{W}_1 \bigl(\bm{x}^{(k)}-\overline{\bm{x}}\bigr), 
\end{equation}
where $\widehat{W}_1\in\R^{(n-r_k)\times n}$ contains the first $n-r_k$ rows of 
$W^{-1}$, and $W_1\in\R^{n\times(n-r_k)}$ is composed of the first $n-r_k$ 
columns of $W$.
\end{theorem}

\smallskip
\begin{proof}
The proof proceeds analogously to that of Theorem 4.2 in~\cite{pr20}.
Replacing $J_k$ and $L$ with their GSVD and setting $\bm{y}=W^{-1}\bm{s}$,
$\bm{z}^{(k)}=W^{-1}\left(\bm{x}^{(k)}-\overline{\bm{x}}\right)$, and
$\bm{g}^{(k)}=U^T\bm{r}_k$,~\eqref{Llinmnls} can be rewritten as 	
the following diagonal least-squares problem
\begin{equation*}%\label{diagprob}
\begin{cases}
\displaystyle \min_{\bm{y}\in\R^n}\|\Sigma_L(\alpha_k\bm{y}+\bm{z}^{(k)})\|^2 \\
\displaystyle \bm{y} \in \bigl\{ \arg \min_{\bm{y}\in\R^n} 
\|\Sigma_J\bm{y}+\bm{g}^{(k)}\|^2 \bigr\}.
\end{cases}
\end{equation*}
When $m\geq n$, the diagonal linear system in the constraint is solved by a 	
vector $\bm{y}$ with entries
$$
y_i = \begin{cases}
\displaystyle
-\frac{g^{(k)}_i}{c_{i-n+r_k\strut}}, \quad & i=n-r_k+1,\ldots,p, \\
-g^{(k)\strut}_i, & i=p+1,\ldots,n.
\end{cases}
$$
The components $y_i$, for $i=1,\ldots,n-r_k$, can be determined by 	
minimizing the norm
\begin{equation}\label{sigmanorm}
\begin{aligned}
\|\Sigma_L(\alpha_k\bm{y}+\bm{z}^{(k)})\|^2 
&= \sum_{i=1}^{n-r_k}\left(\alpha_k y_i+z_i^{(k)} \right)^2 \\
&\phantom{=}
+ \sum_{i=n-r_k+1}^{p} \left( -\alpha_k\frac{g^{(k)}_i}{\gamma_{i-n+r_k}}
+ s_{i-n+r_k} z_i^{(k)} \right)^2,
\end{aligned}
\end{equation}
where $\gamma_i=\frac{c_i}{s_i}$ are the generalized singular
values of the matrix pair $(J_k,L)$.
The minimum of~\eqref{sigmanorm} is reached for 
$
y_i=-\frac{1}{\alpha_k}z^{(k)}_i
=-\frac{1}{\alpha_k}\widehat{\bm{w}}^{i}(\bm{x}^{(k)}-\overline{\bm{x}})
$, $i=1,\ldots,n-r_k$,
and the solution to~\eqref{Llinmnls}, that is, the next approximation to the
solution of~\eqref{Lnonlinmnls}, is 
\begin{equation}
\begin{aligned}
\label{mwngnit}
\bm{x}^{(k+1)} &= \bm{x}^{(k)} + \alpha_k W\bm{y} \\
&= \bm{x}^{(k)} -\sum_{i=1}^{n-r_k} 
z_i^{(k)}\bm{w}_i 
-\alpha_k\sum_{i=n-r_k+1}^{p} \frac{g^{(k)}_i}{c_{i-n+r_k}} \bm{w}_i
-\alpha_k\sum_{i=p+1}^n g^{(k)}_i \bm{w}_i,
\end{aligned}
\end{equation}
where the first summation in the right-hand side can be rewritten as
$W_1\widehat{W}_1(\bm{x}^{(k)}-\overline{\bm{x}})$.
Applying the same procedure to~\eqref{Lmnls},
we obtain
$$
\widetilde{\bm{x}}^{(k+1)} = \bm{x}^{(k)} 
-\alpha_k\sum_{i=n-r_k+1}^{p} \frac{g^{(k)}_i}{c_{i-n+r_k}} \bm{w}_i
-\alpha_k\sum_{i=p+1}^n g^{(k)}_i \bm{w}_i,
$$
from which~\eqref{mlnsol} follows.
Since solving~\eqref{Llinmnls} for $m<n$ leads to a formula similar
to~\eqref{mwngnit}, with $g^{(k)}_{i-n+m}$ in place of $g^{(k)}_i$, the 
validity of~\eqref{mlnsol} is confirmed. $\qquad$
\end{proof}

As in the computation of the minimal-norm solution, the iteration based on
\eqref{mlnsol} fails to converge without a suitable relaxation parameter
$\beta_k$ for the projection vector
$\bm{t}^{(k)}=W_1\widehat{W}_1 (\bm{x}^{(k)}-\overline{\bm{x}})$.
We adopted an iteration similar to \eqref{mngn2}, choosing $\beta_k$ by
adapting Algorithms~\ref{algobeta} and~\ref{algoeta} to this setting.
It is important to note that $\widetilde{\cP}_{\cN(J_k)}=W_1\widehat{W}_1$ is
an oblique projector onto $\cN(J_k)$.

At the same time, the rank of the Jacobian is estimated at each step by
applying the procedure described in Section~\ref{rankjac} to the diagonal
elements $c_j^{(k)}$, $j=1,\ldots,q-d$, of the GSVD factor $\Sigma_J$ of 
$J_k$; see \eqref{csmat}.
In this case, at each step, we compute the ratios
\[
\rho_i^{(k)} = \frac{c_{i+1}^{(k)}}{c_i^{(k)}}, 
\qquad i=1,2,\ldots,q-d-1,
\]
where $q=\min(m,n)$.

Actually, the GSVD routine computes the matrix $W^{-1}$, but the matrix $W$ is
needed for the computation of both the vectors $\widetilde{\bm{s}}^{(k)}$ and
$\bm{t}^{(k)}$.
To reduce the computational load, we compute at each iteration the
LU factorization $PW^{-1}=LU$, and we use it to solve the linear system with
two right-hand sides
$$
W^{-1} \begin{bmatrix} \bm{t}^{(k)} & \widetilde{\bm{s}}^{(k)} \end{bmatrix}
= \begin{bmatrix} 
\widehat{W}_1 (\bm{x}^{(k)}-\overline{\bm{x}}) & \bm{0}_{n-r} \\
\bm{0}_r & \widetilde{\bm{y}}
\end{bmatrix},
$$
where $\widetilde{\bm{y}}\in\R^r$ contains the last $r$ components of the
vector $\bm{y}$ appearing in \eqref{mwngnit}, and $\bm{0}_k$ denotes the zero
vector of size $k$.

\section{Test problems and numerical results}\label{examples}

The MNGN2 method, defined by \eqref{mngn2}, was implemented in the Matlab
programming language; the software is available from the authors.
The developed functions implement all the variants of the MNGN2 algorithm, as
well as the MNGN and CKB methods developed in \cite{pr20} and \cite{campbell},
respectively.

In the following, the MNGN2 algorithm~\eqref{mngn2} will be denoted by different
names, according to the particular implementation.
In the method denoted by MNGN$2_\alpha$, we let $\beta_k=\alpha_k$ in
\eqref{mngn2}, and determine $\alpha_k$ by the Armijo--Goldstein principle.
Algorithm~\ref{algobeta} is denoted by MNGN$2_{\alpha\beta}$, when
$\delta(\rho,\eta)=\eta\rho$, with a fixed value of $\eta$.
The same algorithm with $\delta(\rho,\eta)=\rho^\eta$, and $\eta$ estimated by
Algorithm~\ref{algoeta}, is labeled as MNGN$2_{\alpha\beta\delta}$.
The algorithm~\eqref{camp} developed in~\cite{campbell} is denoted by
CKB$_1$ when $\gamma_k=(0.5)^{k+1}$, and by
CKB$_2$ when $\gamma_k=(0.5)^{2^k}$.
The same algorithms are denoted by rCKB$_1$ and rCKB$_2$ when they are applied
with the automatic estimation of the rank of the Jacobian, discussed in
Section~\ref{rankjac}.
To compare the methods and investigate their performance, we performed numerical
experiments on various test problems that highlight particular difficulties in
the computation of the minimal-norm solution.
Example \ref{exam_robot} illustrates a situation where the MNGN method either
fails or produces unacceptable results, while the other methods perform well;
in Example \ref{exam_camp}, we investigate the dependence of the
MNGN$2_{\alpha\beta}$ method on the choice of the parameter $\eta$;
Example \ref{exam_2} is the first medium-size test problem we consider, it
shows the importance of the Jacobian rank estimation for the effectiveness of
the algorithms;
in Example \ref{exam_1}, the methods are compared in the solution of
minimal-$L$-norm problems with different regularization matrices;
finally, in Example \ref{exam_6}, we let the dimension of the problem vary and
we explore the dependence of the computed solution on the availability of a
priori information in the form of a model profile.

For each experiment, we repeated the computation 100 times, varying the
starting point $\bm{x}^{(0)}$ by letting its components be uniformly
distributed random numbers in $(-5,5)$.
The model profile $\overline{\bm{x}}$ was set to the zero vector except in
Example~\ref{exam_6}.

We consider a numerical test a ``success'' if the algorithm converges according
to condition \eqref{conver}, with stop tolerance $\tau=10^{-8}$ and maximum
number of iterations $N_{\text{max}}=500$.
A failure is not a serious problem, in general, because nonconvergence
simply suggests to try a different starting vector. Anyway, if this happens too
often, it increases the computational load.
At the same time, a success of a method does not imply that it recovers the
minimal-norm solution, as the convergence is only local.
So, to give an idea of the performance of the methods, we measure over all the
tests the average of both the number of iterations required and the norm of
the converged solution $\|\widetilde{\bm{x}}\|$.
We also report the number of successes.

We note that the computational cost of each iteration is roughly the same
for all the methods considered.
Indeed, the additional complexity required by the MNGN2 algorithms consists of
the estimation of the numerical rank $r_{\epsilon,k}$, 
of the residual increase $\delta(\rho,\eta)$, and
of the projection parameter $\beta_k$.
All these computations involve a small number of floating point operations; see
also Remark~\ref{complexity}.
\medskip

\begin{example}\rm\label{exam_robot}
In this first example we consider a nonlinear model that describes the
behavior of a redundant parallel robot. It is a problem that concerns the
inverse kinematics of position, and is defined by the following function
$F:\R^4 \rightarrow \R^2$ 
\begin{equation*}
F(\bm{x})=
\begin{bmatrix}
(X-A\cos(x_1))^2+(Y-A\sin(x_1))^2-x_2^2 \\
(X-A\cos(x_3)-H)^2+(Y-A\sin(x_3))^2-x_4^2
\end{bmatrix},
\end{equation*}
with the data vector $\bm{b}=\bm{0}$ in \eqref{nonlinL2}.
%$$
%\begin{aligned}
%x_1=a\cos(q_1), \quad y_1=a\sin(q_1),\\
%x_2=a\cos(q_3), \quad y_2=a\sin(q_3),
%\end{aligned}
%$$
The model describes the kinematic of a robotic arm moved by 4 motors, whose
position is identified by the unknowns $\{x_i\}_{i=1}^4$, which must reach a
point with given coordinates $(X,Y)$; $A$ and $H$ are parameters describing the
system. In our simulation we assume $(X,Y)=(3,3)$, $A=2$, $H=10$.

The Jacobian matrix of $F$ is
\begin{equation*}
J(\bm{x})= \begin{bmatrix}
\dfrac{\partial F_1}{\partial x_1} & \dfrac{\partial F_1}{\partial x_2} & 0 & 0 
\\
0 & 0 & \dfrac{\partial F_2}{\partial x_3} & \dfrac{\partial F_2}{\partial x_4}
\end{bmatrix},
\end{equation*}
with
$$
\begin{aligned}
\frac{\partial F_1}{\partial x_1} &= 2A(X-A\cos(x_1))\sin(x_1)
	-2A(Y-A\sin(x_1))\cos(x_1),\\
\frac{\partial F_2}{\partial x_3} &= 2A(X-A\cos(x_3)-H)\sin(x_3)
	-2A(Y-A\sin(x_3))\cos(x_3),\\
\frac{\partial F_1}{\partial x_2} &= -2x_2,\quad
\frac{\partial F_2}{\partial x_4} = -2x_4.
\end{aligned}
$$

The results obtained are reported in Table~\ref{tabexrob}.
We see that the MNGN$2_\alpha$ and CKB$_1$ methods recover solutions with
smaller norms, in the average, but the first one requires a large number of
iterations. 
The MNGN$2_{\alpha\beta\delta}$ implementation, with automatic estimation of
the projection step $\beta_k$, quickly converges but produces solutions with
slightly larger norms.
The CKB$_2$ method leads to solutions with a worse norm, testifying that the
performance of the method in \eqref{camp} is very sensitive to the choice of
the sequence $\gamma_k$.
The MNGN method from \cite{pr20} leads to solutions far from optimality, and
fails in 70\% of the tests.
This happens in most of the examples considered in this paper, so we will
involve it only in another experiment.

\begin{table}[ht]\centering
\caption{Results for Example~\ref{exam_robot}.}\label{tabexrob}
\begin{tabular}{lccc}
\hline
method & iterations & $\|\widetilde{\bm{x}}\|$ & \#success \\
\hline
MNGN$2_\alpha$ & 239 & 8.7246 & 92 \\
MNGN$2_{\alpha\beta\delta}$ & 38 & 9.0621 & 96 \\
CKB$_1$ & 26 & 8.5515 & 100 \\
CKB$_2$ & 10 & 9.7344 & 100 \\
MNGN & 182 & 17.6329 & 30 \\
\hline
\end{tabular}
\end{table}

\end{example}
\medskip

\begin{example}\rm\label{exam_camp}
Here we consider a test problem introduced in~\cite{campbell}.
Let $F:\R^3 \rightarrow \R$ be the nonlinear function defined by
\[
F(\bm{x})= x_3 - (x_1-1)^2 - 2(x_2-2)^2 -3.
\]
The equation $F(\bm{x})=0$ represents an elliptic paraboloid in $\R^3$ with 
vertex $\bm{V}=(1,2,3)^T$.
We remark that the minimal-norm solution is the point
$$
\bm{x}^\dagger \approx (0.859754, 1.849178, 3.065164)^T,
$$
and not the vector $\widehat{\bm{x}}$ reported in~\cite[Sec.~4.2]{campbell}. 
Indeed, $\|\bm{x}^\dagger\|\approx 3.681558$, whereas 
$\|\widehat{\bm{x}}\|\approx 3.706359$. 

The results obtained are reported in Table~\ref{tabexcamp}.
The MNGN$2_{\alpha\beta}$ method is tested with two values of the parameter
$\eta$ appearing in the residual increase $\delta(\rho,\eta)=\eta\rho$; see
Algorithm~\ref{algobeta}. It is clear that it can lead to accurate solutions
only if the parameter is suitably chosen ($\eta=2$).
On the contrary ($\eta=8$), it shows a great number of failures.

As in the previous example, the best results are produced by MNGN$2_{\alpha}$,
and MNGN$2_{\alpha\beta\delta}$ reaches very similar solutions but is about 10
times faster.
The CKB methods take a smaller number of iterations, but produce less accurate
solutions.

\begin{table}[ht]\centering
\caption{Results for Example~\ref{exam_camp}.}\label{tabexcamp}
\begin{tabular}{lccc}
\hline
method & iterations & $\|\widetilde{\bm{x}}\|$ & \#success \\
\hline
MNGN$2_{\alpha\beta}\,(\eta=8)$ & 174 & 3.6903 & 15 \\
MNGN$2_{\alpha\beta}\,(\eta=2)$ & 62 & 3.7120 & 100 \\
MNGN$2_\alpha$ & 330 & 3.6816 & 100 \\
MNGN$2_{\alpha\beta\delta}$ & 37 & 3.6832 & 100 \\
CKB$_1$ & 26 & 3.7343 & 100 \\
CKB$_2$ & 10 & 3.7561 & 100 \\
\hline
\end{tabular}
\end{table}

\end{example}
\medskip

\begin{example}\rm\label{exam_2}
Let $F:\R^n \rightarrow \R^m$ be the nonlinear function
\begin{equation}\label{nonlinfun}
F(\bm{x})=\left[ F_1(\bm{x}),F_2(\bm{x}),\ldots,F_m(\bm{x}) \right]^T, \qquad 
m\leq n,
\end{equation}
defined by
\[
F_i(\bm{x}) = \frac{1}{2} S(\bm{x}) \left(x_i^2+1\right),
\qquad i=1,\ldots,m,
\]
where
\[
S(\bm{x}) = \sum_{j=1}^{n} \left(\frac{x_j-c_j}{a_j}\right)^2 -1
\]
is the $n$-ellipsoid with center $\bm{c}=(c_1,\ldots,c_n)^T$ and whose semiaxes 
are the components of the vector $\bm{a}=(a_1,\ldots,a_n)^T$.
The locus of the solutions is the $n$-ellipsoid.

Setting $y_i=x_i^2+1$, for $i=1,\ldots,m$, and $z_j=\frac{x_j-c_j}{a_j^2}$, for 
$j=1,\ldots,n$, the Jacobian matrix can be expressed as
\[
J(\bm{x})=S(\bm{x}) D_{m,n}(\bm{x}) + \bm{y}\bm{z}^T,
\]
where $D_{m,n}(\bm{x})$ is an $m\times n$ diagonal matrix whose main diagonal 
consists of the vector $\bm{x}$. 
Indeed,
\[
\frac{\partial F_i}{\partial x_k}=
\begin{cases}
x_i S(\bm{x}) + \dfrac{x_i-c_i}{a_i^2}\left(x_i^2+1\right), \qquad & 
k=i, \\
\dfrac{x_k-c_k}{a_k^2}\left(x_i^2+1\right), \qquad & k\neq i.
\end{cases}
\]
When $S(\bm{x})=0$, $\rank(J(\bm{x}))=1$, so we expect the Jacobian to be
rank-deficient in a neighborhood of the solution.

If $\bm{a}=\bm{e}=(1,\ldots,1)^T$, the locus of the solutions is the $n$-sphere 
centered in $\bm{c}$ with unitary radius.
If $\bm{c}=2\bm{e}$, the minimal-norm solution is
\[
\bm{x}^\dagger = \left( 2 - \frac{\sqrt{n}}{n}\right) \bm{e},
\]
while if $\bm{c}=(2,0,\ldots,0)^T$ it is $\bm{x}^\dagger=(1,0,\ldots,0)^T$.

Table~\ref{tabrank} displays the results for the last case, when $m=8$ and
$n=10$.
These results aim at underlining the importance of estimating the rank of the
Jacobian $J_k$.
The implementations of the MNGN2 algorithm are more or less equivalent,
recovering solutions with almost optimal norm; MNGN$2_\alpha$ fails in 17\% of
the tests. The value of $\eta$ for MNGN$2_{\alpha\beta}$ is tailored to
maximize the performance, which is not possible in practice, while it is
automatically estimated for MNGN$2_{\alpha\beta\delta}$.
The MNGN and CKB methods do not perform well, because of the rank deficiency of
the Jacobian. We also implemented the rank estimation in the algorithms from
\cite{campbell}; the corresponding methods are denoted by rCKB.
It happens that rCKB$_2$ produces results comparable to the MNGN2
methods, confirming that a correct estimation of the rank is essential for the
convergence, while rCKB$_1$ converges only in 32\% of the tests and 
produces solutions with large norms.
Again, this shows that the sequence adopted for the step length in (r)CKB
methods is critical for the effectiveness of the computation.

\begin{table}[ht]\centering
\caption{Results for Example~\ref{exam_2} with $m=8$, $n=10$, $\bm{a}=\bm{e}$,
and $\bm{c}=(2,0,\ldots,0)^T$. In MNGN, CKB$_1$, and CKB$_2$, the rank is not 
estimated.}\label{tabrank}
\begin{tabular}{lccc}
\hline
method & iterations & $\|\widetilde{\bm{x}}\|$ & \#success \\
\hline
MNGN$2_\alpha$ & 209 & 1.0263 & 83 \\
MNGN$2_{\alpha\beta}\,(\eta=8)$ & 208 & 1.0449 & 99 \\
MNGN$2_{\alpha\beta\delta}$ & 206 & 1.0367 & 97 \\
MNGN & 70 & 2.1083 & 2 \\
CKB$_1$ & 216 & 2.2002 & 32 \\
CKB$_2$ & 20 & 2.1305 & 2 \\
rCKB$_1$ & 160 & 2.1088 & 32 \\
rCKB$_2$ & 197 & 1.0454 & 97 \\
\hline
\end{tabular}
\end{table}

The norms of the solutions, whose average is displayed in Table~\ref{tabrank},
are reported in the boxplot in the left pane of Figure~\ref{ex2boxnorm}.
In each box, the red mark is the median, the edges of the blue box
are the 25th and 75th percentiles, and the black whiskers extend to the most
extreme data points non considered to be outliers, which are plotted as red
crosses.

\begin{figure}
\centering
\includegraphics[width=.47\textwidth]{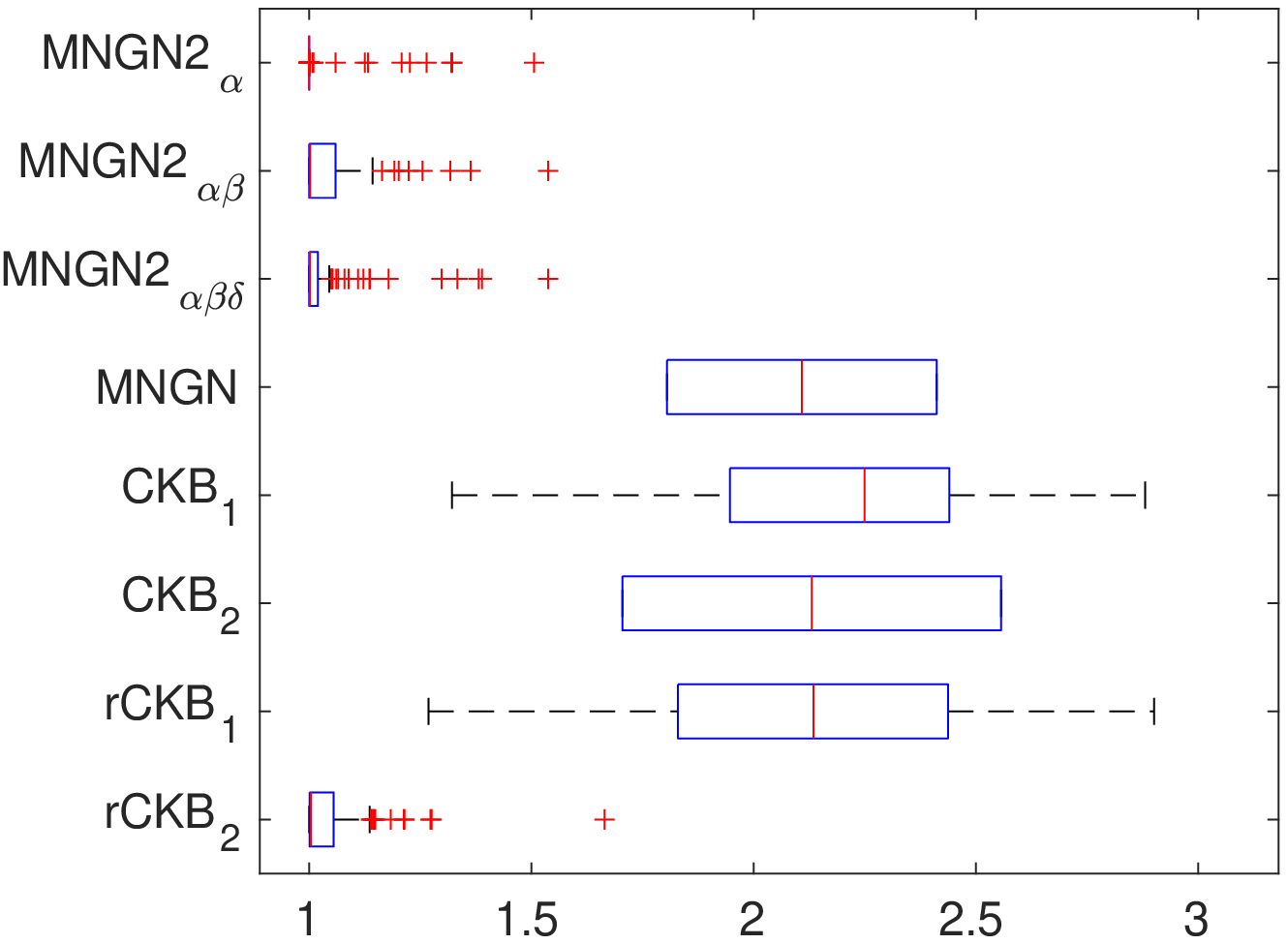} \hspace{0.2cm}
\includegraphics[width=.49\textwidth]{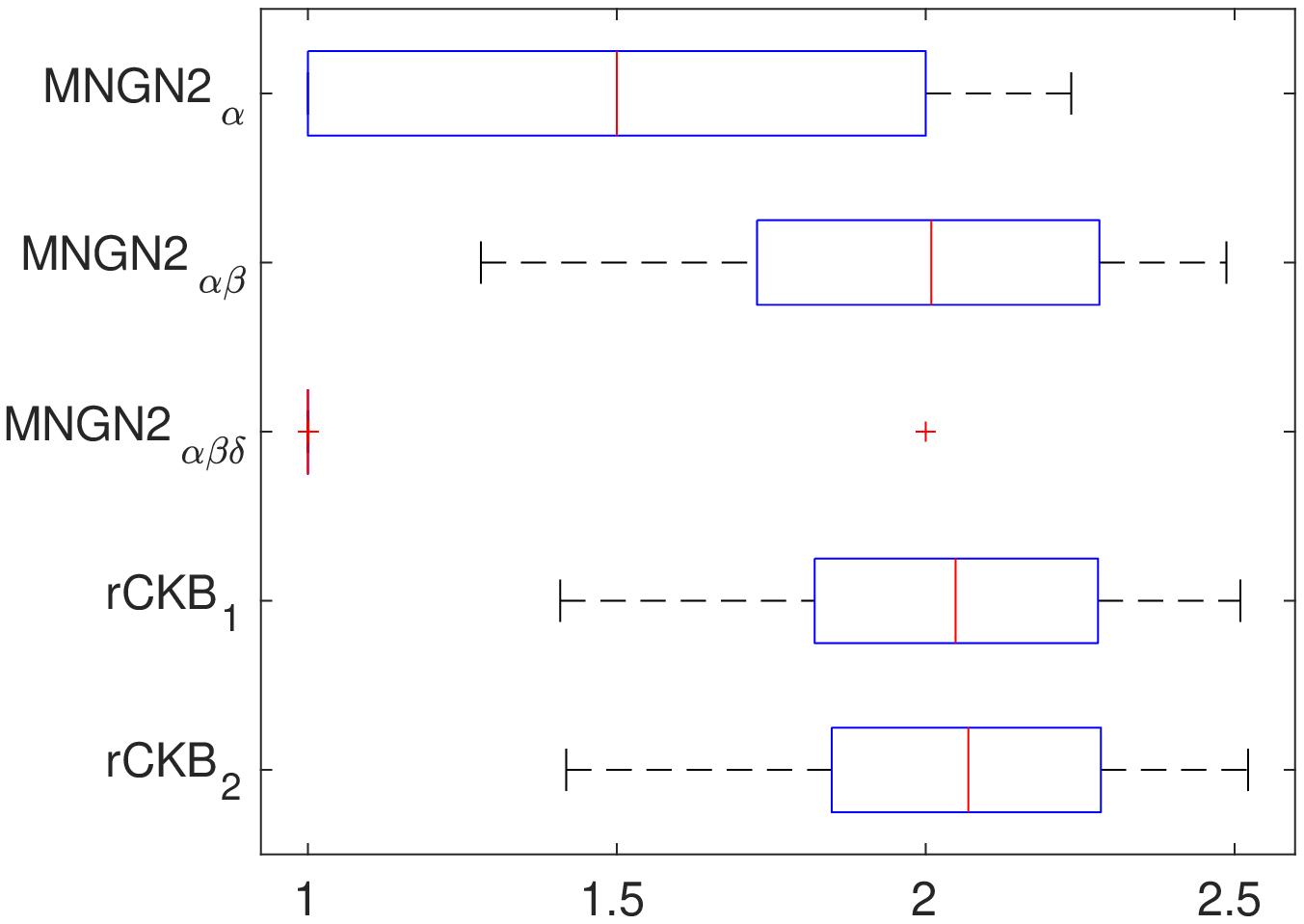} 
\caption{Boxplot of the norms of the solutions for Examples~\ref{exam_2} (left) 
and~\ref{exam_1} (right). The series, labeled by the methods name, are
displayed in the same order of Table~\ref{tabrank} and Table~\ref{tabex1.1},
respectively.}
\label{ex2boxnorm}
\end{figure}

\end{example}
\medskip

\begin{example}\rm\label{exam_1}
Let $F$ be a nonlinear function such as \eqref{nonlinfun}, with
\begin{equation}\label{ex3}
F_i(\bm{x})= S(\bm{x}) \left(x_i-c_i\right), \qquad i=1,\ldots,m,
\end{equation}
and $S(\bm{x})$ defined as in the previous example.
The first order derivatives of $F_i(\bm{x})$ are
\[
\frac{\partial F_i}{\partial x_k}=
\begin{cases}
\dfrac{2}{a_i^2}(x_i-c_i)^2 + S(\bm{x}), \qquad &k=i, \\
\dfrac{2^{\strut}}{a_k^2}(x_k-c_k)(x_i-c_i), \qquad &k\neq i.
\end{cases}
\]
Setting $y_i=x_i-c_i$, for $i=1,\ldots,m$, and $z_j=\frac{x_j-c_j}{a_j^2}$, for 
$j=1,\ldots,n$, the Jacobian matrix can be represented as
\[
J(\bm{x})=S(\bm{x}) I_{m\times n} + 2\bm{y}\bm{z}^T,
\]
where $I_{m\times n}$ includes the first $m$ rows of an identity matrix of size
$n$.
The Jacobian turns out to be a diagonal plus rank-1 matrix.
This structure may be useful to reduce complexity when solving large scale
problems.

When $S(\bm{x})=0$, the matrix $J(\bm{x})$ has rank 1.
Indeed, in this case, the compact SVD of the Jacobian is
\[
J(\bm{x})=\frac{\bm{y}}{\|\bm{y}\|} (2\|\bm{y}\|\|\bm{z}\|) 
\frac{\bm{z}^T}{\|\bm{z}\|},
\]
so that the only non-zero singular value is $2\|\bm{y}\|\|\bm{z}\|$.
As in the preceding example,
we may assume that the Jacobian is rank-deficient in the
surroundings of a solution.
%Therefore, the pseudo-inverse is
%\[
%J(\bm{x})^\dagger= \frac{\bm{z}}{\|\bm{z}\|} 
%\left(\frac{1}{2\|\bm{y}\|\|\bm{z}\|} \right) \frac{\bm{y}^T}{\|\bm{y}\|}.
%\]

\begin{figure}
\centering
\includegraphics[width=\textwidth]{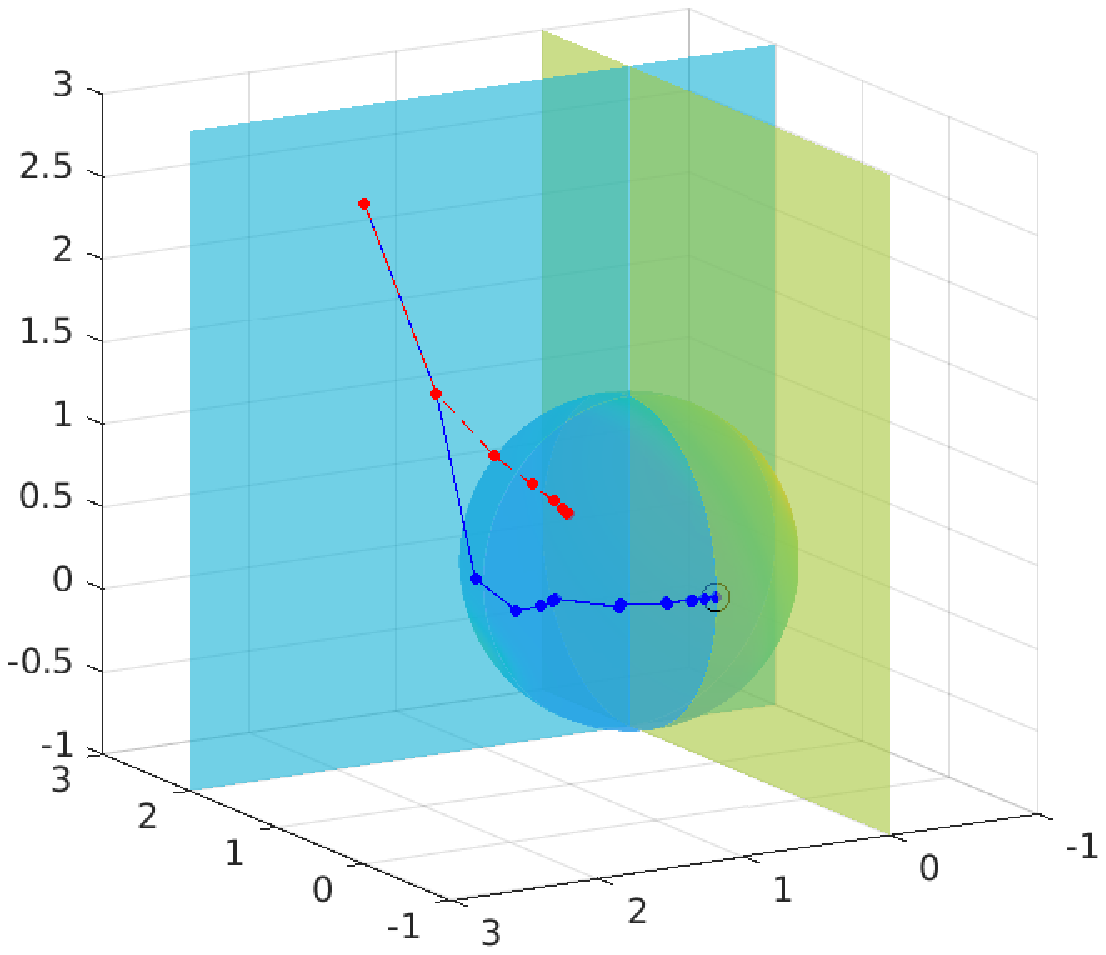} 
\caption{Solution of problem~\eqref{ex3} (Example~\ref{exam_1}) for $m=2$
and $n=3$, with $\bm{a}=(1,1,1)^T$, $\bm{c}=(2,0,0)^T$, and
$\bm{x}^{(0)}=(0,3,3)^T$. The locus of the solutions is the sphere and the line
intersection of the two planes. The blue dots are the iterations of the
MNGN$2_{\alpha\beta\delta}$ method, and the red ones correspond to the rCKB$_1$
method. The black circle encompasses the minimal-norm solution.}
\label{ex3fig}
\end{figure}

The locus of the solutions is the union of the $n$-ellipsoid and the 
intersection between the planes $x_i=c_i$, $i=1,\ldots,m$.

If $\bm{a}=\bm{e}$ and $\bm{c}=2\bm{e}$, the minimal-norm solution 
$\bm{x}^\dagger$ depends on the dimensions $m$ and $n$: if
$m<n-\sqrt{n}+\frac{1}{4}$, then it is
\[
\bm{x}^\dagger=(\underbrace{2,2,\ldots,2}_m,\underbrace{0,\ldots,0}_{n-m})^T,
\]
otherwise, it is
\begin{equation}\label{smoothsol}
\bm{x}^\dagger = \left( 2 - \frac{\sqrt{n}}{n}\right) \bm{e}.
\end{equation}
If $\bm{c}=(2,0,\ldots,0)^T$, it is $\bm{x}^\dagger=(1,0,\ldots,0)^T$.
The case $m=2$, $n=3$, is displayed in Figure~\ref{ex3fig}, together with the
iterations of the algorithms MNGN$2_{\alpha\beta\delta}$ and rCKB$_1$.
In this test, the latter algorithm converges to a solution of non-minimal norm.

Table~\ref{tabex1.1} illustrates the situation where $\bm{a}=\bm{e}$,
$\bm{c}=(2,0,\ldots,0)^T$, $m=8$ and $n=10$. The corresponding boxplot of the 
norms of the solutions is
displayed in the right pane of Figure~\ref{ex2boxnorm}.
The MNGN$2_{\alpha\beta\delta}$ method is the only one which recovers the
correct solution; MNGN$2_{\alpha}$ gets close to it, but with a very small
number of successes.

\begin{table}[ht]\centering
\caption{Results for Example~\ref{exam_1} with $m=8$, $n=10$, $\bm{a}=\bm{e}$,
and $\bm{c}=(2,0,\ldots,0)^T$.}\label{tabex1.1}
\begin{tabular}{lccc}
\hline
method & iterations & $\|\widetilde{\bm{x}}\|$ & \#success \\
\hline
MNGN$2_\alpha$ & 215 & 1.5196 & 12 \\
MNGN$2_{\alpha\beta}\,(\eta=8)$ & 11 & 1.9911 & 100 \\
MNGN$2_{\alpha\beta\delta}$ & 47 & 1.0100 & 100 \\
rCKB$_1$ & 27 & 2.0346 & 100 \\
rCKB$_2$ & 11 & 2.0531 & 100 \\
\hline
\end{tabular}
\end{table}

Table~\ref{tabex1} reports the results obtained for 
$\bm{a}=\bm{e}$ and $\bm{c}=2\bm{e}$.
In this case, the solution is \eqref{smoothsol}.
We applied the algorithms to both the solution of the minimal-norm problem, and
the computation of the minimal-$L$-norm solution with $L=D_2$, i.e., the
discrete approximations of the second derivative \eqref{d1d2}.
Since the solution is exactly in the null space of $L$, we expect the
minimal-$L$-norm solution to perform well.
No algorithm is accurate when $L=I$, as the minimal norm is
$2\sqrt{n}-1=5.3246$.
When $L=D_2$, the two MNGN2 implementations are superior to the rCKB methods, as
$\|L\bm{x}^\dagger\|=0$.
As in the previous example, MNGN$2_{\alpha}$ exhibits a large number of
failures.

\begin{table}[ht]\centering
\caption{Results for Example~\ref{exam_1} with $m=8$, $n=10$, $\bm{a}=\bm{e}$,
and $\bm{c}=2\bm{e}$.}\label{tabex1}
\begin{tabular}{llccc}
\hline
& method & iterations & $\|L\widetilde{\bm{x}}\|$ & \#success \\
\hline
$L=I$
& MNGN$2_\alpha$ & 12 & 5.6569 & 23 \\
%& MNGN$2_{\alpha\beta}$ & 12 & 5.6294 & 100 \\
& MNGN$2_{\alpha\beta\delta}$ & 45 & 5.4529 & 100 \\
& rCKB$_1$ & 26 & 5.7274 & 100 \\
& rCKB$_2$ & 11 & 5.7520 & 100 \\
\hline
$L=D_2$
& MNGN$2_\alpha$ & 20 & 0.0500 & 26 \\
%& MNGN$2_{\alpha\beta}$ & 14 & 0.3956 & 100 \\
& MNGN$2_{\alpha\beta\delta}$ & 17 & 0.0765 & 100 \\
& rCKB$_1$ & 27 & 2.1694 & 100 \\
& rCKB$_2$ & 17 & 2.2761 & 100 \\
\hline
\end{tabular}
\end{table}

Since this example is interesting in itself as a test problem, we report some
further comments on it.
If $m=n$, the locus of the solutions is the union of the $n$-ellipsoid and the 
point $\bm{x}=\bm{c}$.
The spectrum of $J(\bm{x})$ is
\[
\sigma(J(\bm{x})) = \left\{ S(\bm{x})+2\bm{y}^T\bm{z}, S(\bm{x}), \ldots,
S(\bm{x}) \right\},
\]
where the eigenvalue $S(\bm{x})$ has algebraic multiplicity $n-1$.
The Jacobian matrix is invertible if and only if $S(\bm{x})\neq 0$.
If this condition is met, the inverse is obtained by the Sherman--Morrison
formula
\[
J(\bm{x})^{-1}= \frac{1}{S(\bm{x})}I_n - 
\frac{2}{S(\bm{x})(S(\bm{x})+2\bm{z}^T\bm{y})}\bm{y}\bm{z}^T.
\]

%If $m=n-1$ and $\bm{a}=\bm{e}$, the locus of the solutions is the union of the
%$n$-sphere and the intersection line among of the planes $x_j=c_j$,
%$j=1,\ldots,m$.
%If $\bm{c}=2\bm{e}$, the minimal-norm solution is
%\[
%\bm{x}^\dagger = \left( 2 - \frac{\sqrt{n}}{n}\right) \bm{e},
%\]
%while if $\bm{c}=(2,0,\ldots,0)^T$ it is $\bm{x}^\dagger=(1,0,\ldots,0)^T$.

%If $m<n-1$ and $\bm{c}=2\bm{e}$, the minimal-norm solution is
%$\bm{x}^\dagger=(2,0,\ldots,0)^T$,
%while if $\bm{c}=(2,0,\ldots,0)^T$ we have
%$\bm{x}^\dagger=(1,0,\ldots,0)^T$.

\end{example}
\medskip

\begin{example}\rm\label{exam_6}
Let $F$ be the nonlinear function~\eqref{nonlinfun} with components
\begin{equation}\label{ex2}
F_i(\bm{x})= \begin{cases}
S(\bm{x}), \qquad &i=1, \\
x_{i-1}(x_i-c_i), \qquad &i=2,\ldots,m,
\end{cases}
\end{equation}
and $S(\bm{x})$ defined as above.
The first order partial derivatives of $F_i(\bm{x})$ are
\[
\frac{\partial F_i}{\partial x_k}=
\begin{cases}
\dfrac{2}{a_k^2}(x_k-c_k), \quad & i=1, \ k=1,\ldots,n, \\
x_i-c_i, \quad & i=2,\ldots,m, \ k=i-1, \\
x_{i-1}, \quad & i=k=2,\ldots,m, \\
0, & \text{otherwise}.
\end{cases}
\]
Setting $z_j=2\frac{x_j-c_j}{a_j^2}$ and $y_j=x_j-c_j$, for $j=1,\ldots,n$, the 
Jacobian matrix of $F$ is 
\begin{equation}\label{jac2}
J(\bm{x})= \begin{bmatrix}
z_1 & z_2 & z_3 & \cdots & z_{m-1}& z_m & \cdots & z_n \\
y_2 & x_1 & & & & & & \\
& y_3 & x_2 & & & & & \\
& & \ddots & \ddots & & & & \\
& & & \ddots & \ddots & & & \\
& & & & y_m & x_{m-1} & & 
\end{bmatrix}.
\end{equation}
The locus of the solutions is the intersection between the hypersurface
defined by $S(\bm{x})=0$ and by the pairs of planes $x_{i-1}=0$, $x_i-c_i=0$,
$i=2,\ldots,m$.

\begin{figure}
\centering
\includegraphics[width=\textwidth]{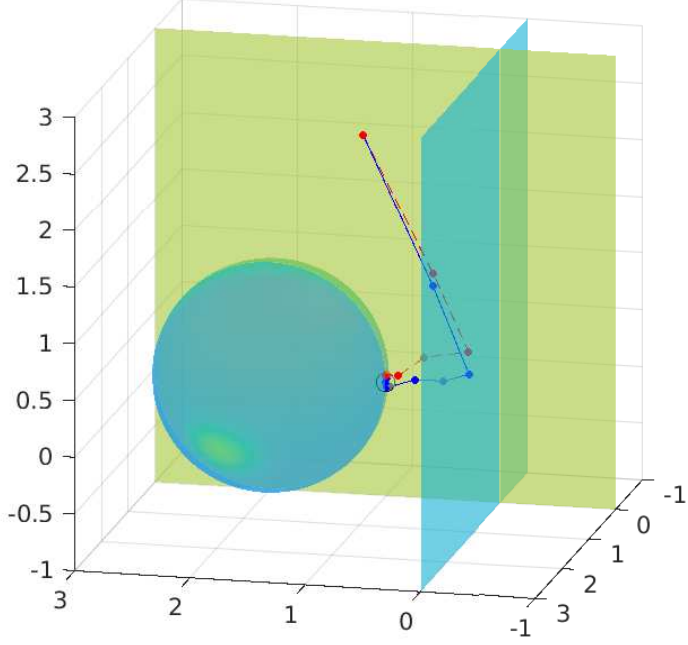}
\caption{Solution of problem~\eqref{ex2} (Example~\ref{exam_6}) for $m=2$ and
$n=3$, with $\bm{a}=(1,1,1)^T$, $\bm{c}=(2,0,0)^T$, and
$\bm{x}^{(0)}=(\frac{1}{2},3,3)^T$. The solutions are in the intersection
between the sphere and the union of the two planes. The blue dots are the
iterations of the MNGN$2_{\alpha\beta\delta}$ method, and the red ones
correspond to the rCKB$_1$ method. The black circle encompasses the minimal-norm
solution.}
\label{ex2fig}
\end{figure}

%If $m=n$, $\bm{a}=\bm{e}=(1,\ldots,1)^T$ and $\bm{c}=2\bm{e}$, the locus 
%of the solutions consists of two points of the $n$-sphere.
%In this case, the minimal-norm solution is $\bm{x}^\dagger=(1,2,\ldots,2)^T$,
%while it is $\bm{x}^\dagger=(1,0,\ldots,0)^T$ if $\bm{c}=(2,0,\ldots,0)^T$.

If $\bm{a}=\bm{e}=(1,\ldots,1)^T$ and $\bm{c}=2\bm{e}$, the minimal-norm 
solution is
%\[
%\bm{x}^\dagger =
%\left( 2-\frac{\sqrt{n-m+1}}{n-m+1}, \underbrace{2, \ldots, 2}_{m-1}, 
%\underbrace{2-\frac{\sqrt{n-m+1}}{n-m+1},\ldots,2-\frac{\sqrt{n-m+1}}{n-m+1}}_{n-m}
%\right)^T,
%\]
\begin{equation}\label{sol1ex2}
\bm{x}^\dagger =
\left( \xi_{n,m}, \underbrace{2, \ldots, 2}_{m-1}, 
\underbrace{\xi_{n,m},\ldots,\xi_{n,m}}_{n-m}
\right)^T,
\end{equation}
with $\xi_{n,m}=2-(n-m+1)^{-1/2}$,
%where the elements equal to $2$ are $m-1$ and the last components 
%$2-\frac{\sqrt{n-m+1}}{n-m+1}$ are $n-m$,
while if $\bm{c}=(2,0,\ldots,0)^T$ it is $\bm{x}^\dagger=(1,0,\ldots,0)^T$.
It is immediate to observe that in the last situation the Jacobian \eqref{jac2}
is rank-deficient at $\bm{x}^\dagger$.
This case is illustrated in Figure~\ref{ex2fig}, where the iterations of
the MNGN$2_{\alpha\beta\delta}$ and the rCKB$_1$ methods are reported too.
The iterations performed are 20 and 24, respectively;
the computed solutions are substantially coincident.

Table~\ref{tabn} displays the results obtained for the same parameter vectors
of Figure~\ref{ex2fig}, when the size of the problem varies, i.e., for
$(m,n)=(8k,10k)$, $k=1,2,3$.
The MNGN2 algorithms behave almost optimally, while the rCKB methods lead to
solutions with larger norm.
The table shows that the performance is not significantly affected by the size
of the problem.
This example suggests that
large scale problems could be faced by the methods discussed, but a suitable
algorithm for the solution of the linearized problem should be adopted, to
reduce the computational complexity of each step.
This aspect will be the object of future research.

\begin{table}[ht]\centering
\caption{Results for Example~\ref{exam_6} with different size $(m,n)$, 
$\bm{a}=\bm{e}$, and $\bm{c}=(2,0,\ldots,0)^T$.}\label{tabn}
\begin{tabular}{llccc}
\hline
$(m,n)$ & method & iterations & $\|\widetilde{\bm{x}}\|$ & \#success \\
\hline
$(8,10)$ & MNGN$2_\alpha$ & 167 & 1.0000 & 48 \\
& MNGN$2_{\alpha\beta}\,(\eta=8)$ & 24 & 1.0508 & 100 \\
& MNGN$2_{\alpha\beta\delta}$ & 37 & 1.0659 & 100 \\
& rCKB$_1$ & 44 & 1.4867 & 100 \\
& rCKB$_2$ & 22 & 1.4776 & 100 \\
\hline
$(16,20)$ & MNGN$2_\alpha$ & 144 & 1.0000 & 36 \\
& MNGN$2_{\alpha\beta}\,(\eta=8)$ & 29 & 1.0170 & 99 \\
& MNGN$2_{\alpha\beta\delta}$ & 34 & 1.0518 & 99 \\
& rCKB$_1$ & 54 & 1.4343 & 100 \\
& rCKB$_2$ & 53 & 1.5269 & 90 \\
\hline
$(24,30)$ & MNGN$2_\alpha$ & 133 & 1.0000 & 34 \\
& MNGN$2_{\alpha\beta}\,(\eta=8)$ & 34 & 1.0154 & 99 \\
& MNGN$2_{\alpha\beta\delta}$ & 32 & 1.0191 & 96 \\
& rCKB$_1$ & 43 & 1.4446 & 100 \\
& rCKB$_2$ & 52 & 1.4529 & 70 \\
\hline
\end{tabular}
\end{table}

Table~\ref{tabxbar} investigates the effectiveness of choosing an appropriate 
model profile $\overline{\bm{x}}$ when applying the MNGN2 algorithms.
We consider the case $\bm{a}=\bm{e}$, $\bm{c}=2\bm{e}$, $m=8$, and $n=10$.
The minimal-norm solution $\bm{x}^\dagger$ is \eqref{sol1ex2}, with
$\xi_{8,10}\simeq 1.4226$ and $\|\bm{x}^\dagger\|\simeq 5.8371$.

When $\overline{\bm{x}}=\bm{0}$, the solutions produced by the considered
variants of the method are almost optimal, but the number of iterations is
quite large, as well as the number of failures for MNGN$2_{\alpha\beta}$ (with a
suitably chosen $\eta$) and MNGN$2_{\alpha\beta\delta}$.
The model profile $\overline{\bm{x}}=2\bm{e}$ reduces the number of iterations
and leads to almost 100\% of successes, but the average norm of the solutions
is slightly larger than the optimal one.
Choosing $\overline{\bm{x}} = 1.7\bm{e}$, a value which is roughly halfway
between 2 and $\xi_{8,10}$, the extreme values of $\bm{x}^\dagger$, restores
the optimality of the results.
This confirms that, when a priori information is available, an accurate choice
of the model profile enhances the performance of the algorithms.

\begin{table}[ht]\centering
\caption{Results for Example~\ref{exam_6} with $m=8$, $n=10$, $\bm{a}=\bm{e}$,
and $\bm{c}=2\bm{e}$.}\label{tabxbar}
\begin{tabular}{llccc}
\hline
& method & iterations & $\|\widetilde{\bm{x}}\|$ & \#success \\
\hline
$\overline{\bm{x}} = \bm{0}$
& MNGN$2_\alpha$ & 138 & 5.8371 & 100 \\
& MNGN$2_{\alpha\beta}\,(\eta=8)$ & 175 & 5.8374 & 38 \\
& MNGN$2_{\alpha\beta\delta}$ & 94 & 5.8988 & 67 \\
\hline
$\overline{\bm{x}} = 2\bm{e}$
& MNGN$2_\alpha$ & 37 & 6.1141 & 99 \\
& MNGN$2_{\alpha\beta}\,(\eta=8)$ & 34 & 6.1144 & 98 \\
& MNGN$2_{\alpha\beta\delta}$ & 34 & 6.1144 & 98 \\
\hline
$\overline{\bm{x}} = 1.7\bm{e}$
& MNGN$2_\alpha$ & 54 & 5.8371 & 100 \\
& MNGN$2_{\alpha\beta}\,(\eta=8)$ & 34 & 5.8394 & 99 \\
& MNGN$2_{\alpha\beta\delta}$ & 40 & 5.8789 & 99 \\
\hline
\end{tabular}
\end{table}

\end{example}

\section{Conclusions}\label{concl}

This paper explores the computation of the minimal-($L$-)norm solution of 
nonlinear least-squares problems, and the reasons for
the occasional lack of convergence of Gauss--Newton methods.
We propose an automatic procedure to estimate the rank of the Jacobian along
the iteration, and the introduction of two different relaxation parameters that
improve the efficiency of the iterative method.
The first parameter is determined by applying the Armijo--Goldstein principle,
while three techniques are investigated to estimate the second one.
In numerical experiments performed on various test problems, the new methods
prove to be very effective, compared to other approaches based on a single
damping parameter.
In particular, the variant which automatically estimates the projection
parameter gives satisfactory results in all the examples.

\section*{Acknowledgements}

The authors are indebted to two anonymous reviewers, whose remarks were
essential for improving both the content and the presentation of this paper.
We thank Maurizio Ruggiu for suggesting the problem reported in
Example~\ref{exam_robot}.
The work of the authors was partially supported by
the Regione Autonoma della Sardegna research project ``Algorithms and 
Models for Imaging Science [AMIS]'' (RASSR57257, intervento finanziato con 
risorse FSC 2014-2020 - Patto per lo Sviluppo della Regione Sardegna),
and the INdAM-GNCS research project ``Tecniche numeriche per l'analisi delle 
reti complesse e lo studio dei problemi inversi''.
Federica Pes gratefully acknowledges CRS4 (Centro di Ricerca, Sviluppo e Studi 
Superiori in Sardegna) for the financial support of her Ph.D. scholarship.

\bibliographystyle{siam}
\bibliography{bibliography}

\end{document}